\documentclass[11pt]{article}

\usepackage[english]{babel}
\usepackage{csquotes}
\usepackage[
    dashed=false,
    style=authoryear,
    maxcitenames=2,
    backend=biber,
    uniquename=false,
    uniquelist=false,
    giveninits=true,
    maxbibnames=99
]{biblatex}

\DeclareNameAlias{sortname}{family-given}
\DeclareLabeldate[misc]{%
  \field{date}
  \field{year}
  \field{eventdate}
  \field{origdate}
  \field{urldate}
}
\renewbibmacro*{volume+number+eid}{%
  \printfield{volume}%
  \iffieldundef{number}
    {}
    {\mkbibparens{\printfield{number}}}%
  \setunit{\addcomma\space}%
  \printfield{eid}}
\renewbibmacro*{in:}{}
\usepackage{amsmath, amsthm}
\usepackage{tikz}
\usepackage{indentfirst}
\usepackage{lipsum}  
\usepackage{amsfonts}
\usepackage{bbm}
\usepackage{hyperref}
\usepackage{graphicx}
\usepackage{subcaption}
\usepackage[toc,page]{appendix}
\usepackage[margin=3.9cm]{geometry}
\usepackage{bbm}
\usepackage{tabularx}
\usepackage{array}
\usepackage{amsmath,amssymb}
\usepackage{mathtools}
\usepackage{enumitem}

\long\def\ignore#1{}

\newtheorem{theorem}{Theorem}
\newtheorem{lemma}{Lemma}

\newtheorem{corollary}{Corollary}

\theoremstyle{definition}
\newtheorem{definition}{Definition}
\newtheorem{assumption}{Assumption}
\newtheorem{proposition}{Proposition}
\theoremstyle{remark}

\newcommand{\cF}{{\cal F}}
\newcommand{\cP}{{\cal P}}

\begin{document} 

\title{\bf Optimal aggregation in multiclass classification}
\author{{\bf Felix Abramovich} \hspace{1cm} {\bf Or Davidovitch-Servadio} \\
Department of Statistics
and Operations Research \\
Tel Aviv University
}
\date{}

\maketitle

\begin{abstract}
This paper studies the problem of aggregation in multiclass classification, generalizing existing results for binary classification. Given a finite dictionary of score functions $\cF=\{f_1,\ldots,f_M\}$ and a training sample, 
the goal is to construct a data-driven aggregate whose risk is as close as possible to that of the best aggregate (oracle) in a prescribed subset of $\operatorname{span}(\cF)$ with respect to a given loss. In particular, we investigate model selection and convex aggregation under the $0/1$ and multinomial logistic losses.

For the $0/1$ loss, we show that the optimal aggregation rates remarkably do not depend on the number of classes $K$, and the corresponding empirical risk minimization (ERM) aggregates are rate-optimal for model selection and nearly-optimal for convex aggregation (up to an additional $\sqrt{\ln(K)}$ factor). 
The situation under the multinomial logistic loss is more subtle. Although the optimal aggregation rates remain independent of $K$, restricting aggregation to proper aggregates in $\operatorname{conv}(\cF)$ leads to a deterioration of the achievable rates with $K$ in certain regimes. 
This reveals a fundamental limitation of proper aggregation in multiclass classification.
To overcome it, we propose improper probability-mixture aggregates that average the induced class probabilities rather than the original score functions and map the result back to score-function space. These aggregates achieve optimal or nearly optimal convergence rates independent of $K$. Finally, we introduce a different aggregation framework in which both the aggregation class and the oracle are defined directly in probability space.

\end{abstract}


\section{Introduction}
Within the aggregation framework, one is given a finite collection of candidate predictors (score functions)
$\cF=\{f_1,\ldots,f_M\}$, referred to as a \emph{dictionary}. 
For a given subset $\Theta \subseteq \mathbb{R}^M$, consider the class of linear combinations $\cF_\Theta=\{f_\theta =\sum_{j=1}^M \theta_j f_j, \; \theta \in \Theta\} \subseteq \operatorname{span}(\cF)$.
The objective is to construct an estimator (aggregate) $\hat{f}$ based on the dictionary and the training sample whose risk is nearly as small as that of the best predictor within $\cF_\Theta$ with respect to the given loss. We do not necessarily require $\hat{f}$ to be a proper aggregate lying in the same class $\cF_\Theta$, nor even in $\operatorname{span}(\cF)$. 
\ignore{
Aggregation can be viewed as a generalization of model selection, where the objective is to select the best performing score function in the dictionary itself.}

More formally,  
let $\mathcal X \times \mathcal Y$ be a measurable space endowed with a probability measure $\pi$, and let
$\mathcal S^n=\{(X_i,Y_i)\}_{i=1}^n,\;
(X_i,Y_i)\overset{\mathrm{i.i.d.}}{\sim}\pi$
be a random sample.
We are given a dictionary of $M$ pre-defined score functions (might be learned from an independent sample) $\cF=\{f_1, ..., f_M\}$, where $f_j:\mathcal X \to \mathcal Z$. Usually,  $\mathcal{X} \subseteq \mathbb{R}^p$ and $\mathcal{Z} \subseteq \mathbb{R}^k$ for some integers $p, k \geq 1$, while $\mathcal{Y}$ is typically $\mathbb{R}^1$
(in a regression setting) or $\{1,\ldots,K\}$ (in a classification setting). 

Let $l:\mathcal Z\times \mathcal Y\to\mathbb R$ be a loss function. The performance of an aggregate $f:\mathcal X\to\mathcal Z$ is measured by its risk $R(f)=\mathbb{E}_{\pi}\big[l(f(X),Y)\big]$ and the goal is to construct $\hat f $ based on the dictionary $\cF$ and the training sample $\mathcal S^n$ satisfying an oracle inequality
\[
    \mathbb{E}_{\pi^{\otimes n}}\big[R(\hat f)\big] 
    \leq 
    \inf_{\theta \in \Theta} R(f_\theta) 
    +
    \epsilon(M,n,\Pi)
\]
with the smallest possible remainder term (excess risk) $\epsilon(M,n,\Pi)$, where the expectation is taken with respect to the sample distribution $\pi^{\otimes n}$ and $\Pi$ is the class of considered distributions. 

The most common choices for $\Theta$ are
\begin{enumerate}
    \item Model selection aggregation: $\Theta$ is the canonical basis $\{e_1,\ldots,e_M\}$ in $\mathbb{R}^M$ or, equivalently,  $\cF_\Theta=\cF$.
    \item Convex aggregation: $\Theta=\Delta_M$, where $\Delta_M$ is the unit simplex, i.e. 
      $\Delta_M=\{\theta \in \mathbb{R}^M: \theta_j \geq 0,\; \sum_{j=1}^{M}\theta_j=1\}$, and, therefore, $\cF_\Theta=\operatorname{conv}(\cF)$, where  $\operatorname{conv}(\cF)=\{f: f=\sum_{j=1}^{M}\theta_jf_j,\;\theta \in \Delta_M\}$ is the convex hull of $
      \cF$.
   \item Linear Aggregation: $\Theta=\mathbb{R}^M$ or, equivalently, $\cF_\Theta=\operatorname{span}(\cF)$.
\end{enumerate}

The ultimate goal is to find an aggregation procedure with the rate-optimal remainder term $\epsilon(M, n, \Pi)$ in the following sense: 
\begin{definition}\label{definition:optimal-rate}
    The remainder term $\epsilon(M, n, \Pi)$ is called  rate-optimal if there exist constants $c_1, c_2 >0$ such that:
    \begin{enumerate}
        \item For any $\cF=\{f_1,...,f_M\}$ there exists an aggregate $\hat{f}$ such that $\forall \pi \in \Pi$, 
        \[
        \mathbb{E}_{\pi^{\otimes n}}[R(\hat{f})]\leq \inf_{\theta \in \Theta}(R(f_\theta)) + c_1\epsilon(M, n,\Pi)
        \]
        \item $\exists \cF=\{f_1,...,f_M\}$ such that for any aggregate $\hat{f}$, $\exists\pi \in \Pi$:
        \[
        \mathbb{E}_{\pi^{\otimes n}}[R(\hat{f})]\geq \inf_{\theta \in \Theta}(R(f_\theta))+c_2\epsilon(M, n, \Pi),
        \]
    \end{enumerate}
\end{definition}

We present below some of the most well-known aggregation procedures.
Let $L^n(f)$ be the empirical loss, i.e.,
\[
L^n(f)=\frac{1}{n}\sum_{i=1}^{n}l(f(X_i),Y_i)
\]

\begin{enumerate}
    \item
    Empirical Risk Minimization (ERM) over $\cF$: 
    \begin{equation} \label{equation:erm-F}
        \hat{f}_{ERM}^{l,\cF}=\operatorname*{arg\,min}_{f \in \cF}L^n(f)
    \end{equation}
    
    \item 
    Empirical Risk Minimization (ERM) over the convex hull of $\cF$: 
    \begin{equation} \label{equation:erm-convex-hull}
        \hat{f}_{ERM}^{l,\operatorname{conv}(\cF)}=\operatorname*{arg\,min}_{f \in \operatorname{conv}(\cF)}L^n(f)
    \end{equation}
    
    \item
    Aggregation with Exponential Weights (AEW): 
    \begin{equation} \label{equation:aew}
        \hat{f}^l_{AEW}=\sum_{j=1}^{M}\hat{\theta}_{AEW, j}f_j, \;\;{\rm where}\;\;\hat{\theta}_{AEW,j}=\frac{e^{-\frac{n}{\beta}L^n(f_j)}}{\sum_{j'=1}^{M}e^{-\frac{n}{\beta}L^n(f_{j'})}}
    \end{equation}
        with a tuning hyperparameter $\beta>0$ called temperature. 

For small $\beta$, the AEW aggregate behaves similarly to the ERM over $\cF$.
The larger $\beta$, the closer the weights are to uniform. 
        
        Note that all vectors in $\Delta_M$ can be viewed as probability distributions over $\cF$. In particular, impose a uniform prior $\pi_j=\frac{1}{M},\;j=1,\ldots,M$ on $\cF$. A straightforward calculation implies that the exponential weights  in \eqref{equation:aew} can be obtained as the solution of the following minimization problem:
        $$
        \hat{\theta}_{AEW}=\operatorname*{arg\,min}_{\theta \in \Delta_M} \left\{\sum_{j=1}^M \theta_j L^n(f_j)+
        \frac{\beta}{n}KL(\theta \parallel \pi)\right\},
        $$
        where the Kullback-Leibler divergence $KL(\theta \parallel \pi)=\sum_{j=1}^M \theta_j \ln\Big(\frac{\theta_j}{\pi_j}\Big)=\sum_{j=1}^M\theta_j \ln\theta_j+\ln M$.
            
    \item Cumulative AEW (CAEW):
    \begin{equation} \label{equation:caew} 
    \hat{f}^l_{CAEW}=\sum_{j=1}^{M}\hat{\theta}_{CAEW,j}f_j, \;\;{\rm where}\;\;\hat{\theta}_{CAEW,j}=\frac{1}{n}\sum_{i=0}^{n-1}\frac{e^{-\frac{i}{\beta}L^{i}(f_j)}}{\sum_{j'=1}^{M}e^{-\frac{i}{\beta}L^{i}(f_{j'})}}
    \end{equation}
    with $L^0(f)=0$. Evidently, $\hat{f}^l_{CAEW}=\frac{1}{n}\sum_{t=i}^{n-1} \hat{f}^l_{AEW,i}$~, where $\hat{f}_{AEW,t}$ is the AEW aggregate \eqref{equation:aew} based on the sample ${\cal S}^t=\{(X_i,Y_i)\}_{i=1}^t$ of size $t$.
    \ignore{
    \textcite{Juditsky_2008} showed the relationship between the weights $\hat{\theta}_{CAEW}$ in \eqref{equation:caew} and mirror averaging in stochastic gradient descent.}

    \item \label{def:q-aggregation} Q-Aggregation:
    \begin{equation} \label{equation:q-aggregation}
        \hat{\theta} = \operatorname*{arg\,min}_{\theta \in \Delta_M}\left\{\nu L^n\left(\sum_{j=1}^{M}\theta_jf_j\right)+(1-\nu)\sum_{j=1}^{M}\theta_jL^n(f_j)+\frac{\beta}{n}\sum_{j=1}^M \theta_j \ln\theta_j\right\},\;\;\; 0 \leq \nu \leq 1,
    \end{equation}
    where the last term, up to an additive constant, is the Kullback-Leibler divergence $KL(\theta \parallel \pi)$ for the uniform distribution $\pi_j=\frac{1}{M},\;j=1,\ldots,M$ (see above). In fact, AEW can be viewed as a particular case of Q-aggregation for $\nu=0$. 
\end{enumerate}

The aggregation problem has been extensively studied in the literature. Early works by \textcite{Nemirovsky_book} and \textcite{Juditsky_and_Nemirovsky} considered aggregation procedures for Gaussian regression under the quadratic loss $l(f(x),y)=(f(x)-y)^2$ and established a convergence rate of $\mathcal{O}\left(\sqrt{\frac{\ln M}{n}}\right)$ for convex aggregation. \textcite{Tsybakov_optimal_rates_regression} extended these results to other types of aggregation and showed that the optimal rates are of the order
$\frac{\ln M}{n}$,
$\min\left\{\frac{M}{n},\sqrt{\frac{\ln\left(1+\frac{M}{\sqrt n}\right)}{n}}\right\}$,
and $\frac{M}{n}$ for model selection, convex, and linear aggregation, respectively. \textcite{Lecu__2013} showed that the aggregate $\hat{f}_{ERM}^{l,\operatorname{conv}(\cF)}$ in \eqref{equation:erm-convex-hull} achieves the optimal rate for convex aggregation, while \textcite{Dai_2012} proved that the Q-aggregate \eqref{equation:q-aggregation} attains the optimal rate for model selection aggregation.

\textcite{Rigollet_2012} extended the results from Gaussian regression to generalized linear models (GLMs) with fixed design under the Kullback--Leibler loss for model selection, convex, and linear aggregation. \textcite{Juditsky_2008} considered a general $\frac{1}{\beta}$-exp-concave loss and showed that, for the CAEW aggregate with temperature $\beta$, the excess risk in model selection aggregation is bounded by $\frac{\beta \ln(M)}{n}$ which is the optimal rate \parencite{lecue2007suboptimalitypenalizedempiricalrisk}. \textcite{Lecue_AEW_low_temp} proved that, under an additional Bernstein condition on the loss, AEW is optimal for model selection aggregation in the low-temperature regime with $\beta<c_1$ for some $c_1>0$. More recently, \textcite{hogsgaard2026aggregationexponentialweightsoptimal} extended these results to higher-temperature regimes and showed that, if the loss is bounded by $B$, $L$-Lipschitz, and $\mu$-strongly convex, the excess risk of the AEW estimator with a constant temperature parameter $\beta>c(B,L,\mu)$ is at most $\frac{\beta\ln(M)}{n}$. On the other hand, they also established the suboptimality of AEW in the high-temperature regime where $\beta_n\to\infty$ as $n$ increases. Since a Lipschitz and strongly convex loss is also exponentially concave \parencite[see, e.g.,  Proposition 3.2 in][] {puchkin2023exploringlocalnormsexpconcave}, these results provide a connection between the optimality properties of AEW and CAEW aggregation. \textcite{lecuerigollet14Qaggregation} generalized the results of \textcite{Dai_2012} for Q-aggregation to general Lipschitz and strongly convex losses when $Y$ is bounded.

In this work, we study aggregation within the classification setting. To the best of our knowledge, all previous results have been obtained for binary classification.
\textcite{lecue2007suboptimalitypenalizedempiricalrisk} studied model selection aggregation for binary classification. He showed that the ERM aggregate \eqref{equation:erm-F} attains the optimal rate \(\sqrt{\ln(M)/n}\) with respect to the $0/1$-loss, and that the ERM \eqref{equation:erm-F}, AEW \eqref{equation:aew}, and CAEW \eqref{equation:caew} aggregates are rate-optimal with respect to the hinge loss. \textcite{lecue_hinge_loss} established the optimality of the ERM \eqref{equation:erm-convex-hull}, AEW \eqref{equation:aew}, and CAEW \eqref{equation:caew} aggregates for convex aggregation under the hinge loss. Both papers also extended these results under the additional margin (low-noise) condition of \textcite{Tsybakov_optimal_rates_classification}. Finally, the general theory of Q-aggregation developed by \textcite{lecuerigollet14Qaggregation} implies that the Q-aggregation procedure is rate-optimal for model selection aggregation under the logistic loss.

In this paper, we generalize aggregation results from binary to multiclass classification. We consider the $0/1$ and multinomial logistic (Kullback--Leibler) losses for model selection and convex aggregation.

For the $0/1$-loss, we show that the optimal aggregation rates do not depend on the number of classes and coincide with those for binary classification. In particular, the ERM aggregates \eqref{equation:erm-F} and \eqref{equation:erm-convex-hull} are rate-optimal respectively for model selection and convex aggregation, up to an additional $\sqrt{\ln(K)}$ factor in the latter case. 

For the multinomial logistic loss, the picture is more subtle. Although the minimax lower bounds for both model selection and convex aggregation are independent of the number of classes $K$, the rates attained by standard proper aggregation procedures may deteriorate with $K$. More importantly, this deterioration is not merely a shortcoming of the particular procedures considered, but a fundamental limitation of proper aggregates in certain regimes. Thus, allowing improper aggregates is essential for attaining the $K$-independent optimal rates in these regimes.

To overcome this limitation, we introduce probability-mixture aggregates that average the class probabilities induced by the original score functions and then map the result back to score-function space. For model selection aggregation, the proposed CAEW probability-mixture aggregate achieves the optimal rate $\frac{\ln(M)}{n}$, independently of $K$. For convex aggregation, a continuous version of this procedure, combined with ERM, yields nearly optimal rates that are also independent of $K$. These results demonstrate that the dependence on $K$ can arise from restricting aggregation to the original score-function space, rather than from the intrinsic statistical complexity of multiclass aggregation.

Finally, we consider a different aggregation framework in which both the aggregation class and the oracle are defined directly in probability space.

The paper is organized as follows.
 Section \ref{section:notations} establishes some basic properties of the multinomial logistic loss that will be used to prove the lower and upper bounds. In Section \ref{section:ms-aggregation}, we study model selection aggregation, derive lower bounds, and present rate-optimal aggregates. Convex aggregation is discussed in Section \ref{section:convex}. A new aggregation framework in probability space is introduced in Section \ref{section:probaggregation}.  All proofs are given in the Appendix.

\section{Notation and Preliminaries} \label{section:notations}

\subsection{Definitions}
Given a measurable  space $\mathcal X \times \{1,\ldots,K\}$ endowed with a probability measure $\pi$, denote by $\pi_{X}$ the marginal distribution on $\mathcal X$. 
We have a dictionary of functions $\cF=\{f_1, ..., f_M\}, \;\;f_j:\mathcal X \to \mathbb{R}^K$. Assume that all $f_j$ are bounded in sup-norm, i.e.
$\lVert f_j(x) \rVert _\infty \leq B$ for any $x \in \mathcal X$, where $B \geq 1$. 

Consider a class of score functions $\cF_\Theta=\{f_\theta=\sum_{j=1}^M \theta_j f_j, \theta \in \Theta\}$, where $\Theta=\{e_1,\ldots,e_M\}$, the canonical basis of $\mathbb{R}^M$, for model selection aggregation, and $\Theta=\Delta_M$ for convex aggregation. 

Within the multiclass framework we consider two losses defined for all functions $f:\mathcal{X} \to \mathbb{R}^K$:
\begin{enumerate}
   \item $0/1$ loss:
        \[
            l_{0/1}(f(x),y) = \mathbb{I}(\operatorname*{arg\,max}_{1 \leq k \leq K} f_k(x) \neq y)
        \]
    \item multinomial logistic loss:
        \begin{gather} \label{equation:multinomial_loss}
            l(f(x),y) = -f_{y}(x)+\ln\left(\sum_{k=1}^{K}e^{f_k(x)}\right) 
        \end{gather}
        where $f_{y}(x)$ is the $y$'th coordinate of $f(x)$.
\end{enumerate}
The multinomial logistic loss is evidently equivalent to the Kullback-Leibler loss $KL(p_y \parallel p_f(x))$ between the empirical distribution of the data, where $p_{yk}=\delta_{yk}$, and the distribution $p_f(x)$ induced by $f(x)$, where $p_f(x)\coloneqq\operatorname{softmax}(f(x))$, i.e., $p_f(x)_k=\frac{e^{f_k(x)}}{\left(\sum_{k'=1}^{K}e^{f_{k'}(x)}\right)}$.

The multinomial logistic loss \eqref{equation:multinomial_loss} is not identifiable without an extra constraint on the functions $f_j$ since it is invariant under shifting all entries $f_{jk}$ by the same arbitrary constant $c$: $l(f(x),y)=l(f(x)+c{\bf 1}_K,y)$.
\ignore{
\[
    \frac{e^{f^{y}(x)+c}}{\sum_{k=1}^Ke^{f^{k}(x)+c}}=\frac{e^{f^{y}(x)}}{\sum_{k=1}^Ke^{f^{k}(x)}}
\]
}

We adopt a symmetric constraint
\begin{equation} \label{equation:symmetric_constraint}
\sum_{k=1}^K f_{jk}(x)=0
\end{equation}
for all $j=1,\ldots,M$ and $x \in \mathcal{X}$, which implies $\sum_{k=1}^K f_{\theta, k}(x)=0$ for all $\theta \in \Theta$.

The corresponding empirical and expected loss functions are
\begin{gather*}
    L_{0/1}^n(f)=\frac{1}{n}\sum_{i=1}^{n}l_{0/1}(f(X_i),Y_i),\; \;\;  R_{0/1}(f)=\mathbb{E}_{(X,Y) \sim \pi}[l_{0/1}(f(X),Y)] = \mathbb{P}_{(X,Y) \sim \pi}(\operatorname*{arg\,max}_{1 \leq k \leq K} f_k(X) \neq Y) \\
    L^n(f)=\frac{1}{n}\sum_{i=1}^{n}l(f(X_i),Y_i), \;\;\;  R(f)=\mathbb{E}_{(X,Y) \sim \pi}[l(f(X),Y)]  
\end{gather*}

Recall that for a given loss $l(f(x),y)$ we want
to construct an aggregated classifier $\hat f $ based on the dictionary of classifiers  $\cF$ and the training sample $\mathcal S^n$ satisfying the (sharp) oracle inequality
\[
    \mathbb{E}_{\pi^{\otimes n}}\big[R(\hat f)\big] 
    \leq 
    \inf_{\theta \in \Theta} R(f_\theta) 
    +
    \epsilon(M,n,\Pi)
\]
with the minimal (rate-optimal) remainder term $\epsilon(M,n,\Pi)$.

Finally, since the AEW and CAEW aggregates are defined throughout the paper using the multinomial logistic loss, we henceforth omit the superscript $l$ and write $\hat{f}_{AEW}$ and $\hat{f}_{CAEW}$ for the corresponding aggregates in  \eqref{equation:aew} and \eqref{equation:caew}, respectively.

\ignore{
\begin{assumption}\label{assumption:identifiability}
    \[
    \forall x \in X, 1\leq j\leq M,\;\; \sum_{k=1}^{K}f_j^k(x)=0 
    \]
\end{assumption}

\begin{assumption}\label{assumption:boundedness}
    \[
    \forall 1\leq j\leq M,\;\; \lVert f_j \rVert_{\infty}\leq1
    \]
\end{assumption}

The first assumption is made for identifiability, since changing all entries by a constant does not change the losses we consider:
\[
    \frac{e^{f^s(x)+c}}{\sum_{k=1}^Ke^{f^k(x)+c}}=\frac{e^ce^{f^s(x)}}{\sum_{k=1}^Ke^ce^{f^k(x)}}=\frac{e^{f^s(x)}}{\sum_{k=1}^Ke^{f^k(x)}}
\]
The boundedness assumption is a natural assumption made in most research in this field.
}

\ignore{
It is easy to see that if the marginal distribution on $\mathcal X$ has a density function, the logistic loss is equivalent to the $KL-Loss$. Take $f^*:\mathcal{X} \to R^K$ that satisfies for every $x,y$: $\pi_{\mathcal Y|\mathcal X=x}(y)=\frac{e^{(f^{*})^{y}(x)}}{\sum_{k=1}^{K}e^{(f^{*})^{k}(x)}}=softmax(f^{*}(x))^y$. For every $f:\mathcal{X} \to R^K$, denote by $\pi_f(x,y)=\pi_X(x)\frac{e^{f^y(x)}}{\sum_{k=1}^{K}e^{f^k(x)}}$. Using this notation we can say that for every $f,g:\mathcal{X} \to R^K$:
\begin{gather*}
    KL(\pi_{f^*} \parallel \pi_{f})-KL(\pi_{f^*} \parallel\pi_{g})=\mathbb{E}_{\pi_{f^*}}[\ln\frac{\pi_{f^*}(x,y)}{\pi_{f}(x,y)}] - \mathbb{E}_{\pi_{f^*}}[\ln\frac{\pi_{f^*}(x,y)}{\pi_{g}(x,y)}] =  \\
    \mathbb{E}_{\pi_{f^*}}[\ln(\pi_{g}(x,y))] - \mathbb{E}_{\pi_{f^*}}[\ln(\pi_{f}(x,y))] = \mathbb{E}_{\pi_{f^*}}[\ln\frac{e^{g^y(x)}}{\sum_{k=1}^{K}e^{g^k(x)}}] - \mathbb{E}_{\pi_{f^*}}[\ln\frac{e^{f^y(x)}}{\sum_{k=1}^{K}e^{f^k(x)}}] = \\
    L(f) -L(g) 
\end{gather*}
}

\subsection{Geometry of the logistic loss}
We now present several important properties of the multinomial logistic loss that will be used to establish risk bounds for the resulting aggregated classifiers.

\begin{definition} \label{definition:matrix_seminorm}
  Let $H \in \mathbb{R}^{m \times m}$ be a symmetric positive semidefinite matrix. The seminorm of $w \in \mathbb{R}^m$ induced by $H$ is $\lVert w \rVert _H = \sqrt{w^T H w}$.
\end{definition}  

\begin{definition}\label{definition:strong-convexity-lipschitz}
    Let $H \in \mathbb{R}^{m \times m}$ be a symmetric positive semidefinite matrix and $A \subseteq \mathbb{R}^m$ a convex set. We say that a twice differentiable function $f:A \to \mathbb{R}$ is 
    \begin{enumerate}
        \item $\lambda$-strongly convex with respect to the  seminorm induced by $H$ if for every $x,\tilde{x} \in A$:
            \begin{gather*}
                f(x)-f(\tilde{x} ) \geq \nabla f(\tilde{x})^T(x-\tilde{x}) +\frac{\lambda}{2}\lVert x-\tilde{x}\rVert_H^2 
            \end{gather*}
        \item $\beta$-smooth with respect to the seminorm induced by $H$ if for every $x,\tilde{x}  \in A$ :
            \begin{gather*}
                f(x)-f(\tilde{x} ) \leq \nabla f(\tilde{x})^T(x-\tilde{x}) +\frac{\beta}{2}\lVert x-\tilde{x}\rVert_H^2 
            \end{gather*}
        \item $L$-Lipschitz with respect to the seminorm induced by $H$ if for every $x,\tilde{x}  \in A$:
             \begin{gather*}
                 |f(x)-f(\tilde{x})| \leq L \lVert x-\tilde{x} \rVert_H
             \end{gather*}
    
    \end{enumerate}
\end{definition}

\begin{definition}\label{definition:exp-concavity} 
    Let $A \subseteq \mathbb{R}^m$ be a convex set. We say that a function $T(x): A \to \mathbb{R}$ is $\beta$-exp-concave if the function $e^{-\beta T(x)}:A\subseteq \mathbb{R}^m \to \mathbb{R}$  is concave.
\end{definition}

\ignore{
    The following lemmas will be very useful in establishing lower and upper bounds for the convergence rates. Since most techniques such as Fano's and Assaoud's require a distance (or at least a semi-distance) over the space of estimators, which the multinomial logistic loss is not, we need to bound it using the $L_2$ metric. First if we consider the loss with respect to a vector function $f(x)$:
}

\begin{lemma}\label{lemma:loss-strongly-convex-and-smooth-wrt-function} Consider the multinomial logistic loss $l(f(x),y)$. For every $x \in \mathcal{X}, y \in \{1,\ldots,K\}$, and $ f_1, f_2:\mathcal{X} \to [-B,B]^K$ satisfying the symmetric constraint \eqref{equation:symmetric_constraint}: 
    \begin{gather*}
        l(f_1(x),y )-l(f_2(x),y) \geq \nabla l(f_2(x),y)^T(f_1(x)-f_2(x))  + \frac{e^{-2B}}{2K}\lVert f_1(x)-f_2(x) \rVert_2^2 \\ 
        l(f_1(x),y)-l(f_2(x),y) \leq \nabla l(f_2(x),y)^T(f_1(x)-f_2(x)) + \frac{e^{2B}}{2K}\lVert f_1(x)-f_2(x)\rVert_2^2
    \end{gather*}
        where the gradient $\nabla l(\cdot,\cdot)$ is taken with respect to the first argument. 
\end{lemma}

Definition \ref{definition:strong-convexity-lipschitz} and Lemma \ref{lemma:loss-strongly-convex-and-smooth-wrt-function} immediately imply the following result:
\begin{lemma}\label{lemma:logistic-loss-geometry}
    For a given dictionary $\cF=\{f_1, \dots f_M\}$ define the matrix $F(x)=(f_1(x),\ldots,f_M(x)) \in \mathbb{R}^{K \times M}$. Then, the multinomial logistic loss $l(F(x)\theta,y)$ in \eqref{equation:multinomial_loss} is $\frac{e^{-2B}}{K}$-strongly convex, $\frac{e^{2B}}{K}$-smooth, and $\sqrt{2}$-Lipschitz over $\Delta_M$ with respect to $\theta$ and the seminorm induced by the symmetric positive semidefinite matrix $H(x)=F(x)^TF(x) \in \mathbb{R}^{M \times M}$, that is, $H_{jl}(x)=f_j(x)^T f_l(x)$.
   Furthermore, it is $\frac{1}{2e^{2B}K}$-exp-concave on $\Delta_M$, that is, the function $e^{-\frac{1}{2e^{2B}K}l(F(x)\theta,y)}$ is concave in $\theta$ on $\Delta_M$.
\end{lemma}

\section{Model selection aggregation} \label{section:ms-aggregation}

In this section, we explore model selection  aggregation, where $\Theta=\{e_1,\ldots,e_M\}$ and, therefore, $\cF_\Theta=\cF$. Recall the score functions are assumed to be uniformly bounded: $\lVert f_j(x) \rVert_{\infty} \leq B$, and to satisfy the symmetric identifiability constraint \eqref{equation:symmetric_constraint}. For a given loss (multinomial logistic or $0/1$) we seek an aggregated classifier $\hat{f}$ based on the training sample $\mathcal S^n$ satisfying
$$    
  \mathbb{E}_{\pi^{\otimes n}}\big[R(\hat f)\big] 
    \leq 
    \min_{1 \leq j \leq M} R(f_j) 
    +
    \epsilon(M,n,\Pi)
$$
with the smallest possible (rate-optimal) remainder term $\epsilon(M,n,\Pi)$.

We start with the ERM aggregate \eqref{equation:erm-F} and establish its upper and lower bounds with respect to the $0/1$ and multinomial logistic losses. We show that while $\hat{f}_{ERM}^{l_{0/1},\cF}$ is rate-optimal with respect to the $0/1$ loss, $\hat{f}_{ERM}^{l,\cF}$ is only suboptimal with respect to the multinomial logistic loss. We then consider the CAEW aggregate \eqref{equation:caew} and show that with a properly chosen tuning parameter $\beta$ it outperforms $\hat{f}_{ERM}^{l,\cF}$ but is still sub-optimal. Moreover, we prove that {\em any proper} aggregate within $\operatorname{conv}(\cF)$ is suboptimal and present an {\em improper} aggregate outside $\operatorname{span}(\cF)$ that does achieve the optimal rate with respect to the multinomial logistic loss. 

\begin{theorem}\label{theorem:erm-upper-bound-expectation}
    Let $\hat{f}_{ERM}^{l_{0/1},\cF}$ and $\hat{f}_{ERM}^{l,\cF}$ be the ERM aggregates \eqref{equation:erm-F} over a dictionary $\cF$ with respect to the $0/1$ loss and multinomial logistic loss, respectively. For every $K,M,n,$ and distribution $\pi$ we have
    \begin{gather} \label{equation:MS01upper}
        \mathbb{E}_{\pi^{\otimes n}}[R_{0/1}(\hat{f}_{ERM}^{l_{0/1},\cF})] \leq \min_{1 \leq j \leq M}R_{0/1}(f_j) + \min\left\{1,4\sqrt{\frac{\ln(M)}{2n}}\right\} 
    \end{gather}
    \begin{gather} \label{equation:MSlogisticupper}
        \mathbb{E}_{\pi^{\otimes n}}[R(\hat{f}_{ERM}^{l,\cF})] \leq \min_{1 \leq j \leq M}R(f_j) + 16B\sqrt{\frac{\ln(M)}{2n}} 
    \end{gather}
    
\end{theorem}
Theorem 2 of \textcite{lecue2007suboptimalitypenalizedempiricalrisk} establishes that for binary classification the optimal rate for model selection aggregation with respect to $0/1$ loss is $\sqrt{\frac{\ln(M)}{n}}$.
Evidently, it cannot be faster for the harder multiclass case. That yields an immediate corollary: 
\begin{corollary}\label{corollary:erm-optimal-0/1}
    The ERM aggregate $\hat{f}_{ERM}^{l_{0/1},\cF}$  is rate-optimal for MS aggregation with respect to the $0/1$ loss, where the optimal rate is $\sqrt{\frac{\ln(M)}{n}}$ regardless of the number of classes $K$.
\end{corollary}

Regarding the AEW and CAEW aggregates,   Corollary 1 of
\textcite{lecue_hinge_loss} implies that even in binary classification, they satisfy a {\em non-sharp} oracle inequality of the form
\[
     \mathbb{E}_{\pi^{\otimes n}}[R_{0/1}(\hat{f}_{AEW/CAEW})] 
     \leq  
     2(1+a)\min_{1 \leq j \leq M} R_{0/1}(f_j) 
     + 
     c(a)\min\left\{1,\sqrt{\frac{\ln(M)}{n}}\right\},
\]
for any $a>0$.

We now discuss the optimality of ERM aggregation with respect to multinomial logistic loss. The following theorem shows that the upper bound \eqref{equation:MSlogisticupper} for $\hat{f}_{ERM}^{l,\cF}$ in Theorem \ref{theorem:erm-upper-bound-expectation} is, in fact, tight and cannot be improved:
\begin{theorem}\label{theorem:erm-lower-bound}
    Assume that $\mathcal{X}$ contains at least $2^M$ distinct points. There exists an integer $M_0 > 0$ such that, for every $K$, $M \geq M_0$, and $n \geq 32 \pi M\ln(M)$, there exists a dictionary $\cF$ of size $M$ and a distribution $\pi$ such that 
    \begin{gather*}
         \mathbb{E}_{\pi^{\otimes n}}
        \Big[ 
            R(\hat{f}_{ERM}^{l,\cF})
            -
            \min_{1 \leq j \leq M}R(f_j)
        \Big] \geq   
        c\sqrt{\frac{\ln(M)}{n}}
    \end{gather*}
    for some $c>0$.
\end{theorem}

The rate $\sqrt{\frac{\ln(M)}{n}}$ attained by ERM aggregation under the logistic loss is, however, not optimal.
Even for the binary case, the optimal aggregation rate is known to be $\frac{\ln(M)}{n}$ \parencite{lecue2007suboptimalitypenalizedempiricalrisk, Rigollet_2012}. 
The inherent limitation of the vanilla ERM aggregation over $\cF$ lies in restricting aggregation to the finite dictionary $\cF$, thereby failing to exploit the convexity of the logistic loss. A natural idea is to enlarge the aggregation class from $\cF$ to its convex hull $\operatorname{conv}(\cF)$ and to consider the corresponding ERM aggregate  $\hat{f}_{ERM}^{l,\operatorname{conv}(\cF)}$ defined in \eqref{equation:erm-convex-hull}.   However, \textcite{Aggregation-via-erm-lecue} showed that the convergence rate for such aggregation is lower bounded by $\sqrt{\frac{1}{n}}$ in the regression setting, which is still substantially slower than the optimal aggregation rate $\frac{\ln(M)}{n}$, and indicates a limitation of this strategy. Passing from $\cF$ to $\operatorname{conv}(\cF)$ improves the geometric structure of the aggregation class by making it convex, but simultaneously significantly increases its complexity. This trade-off suggests that achieving the optimal aggregation rate requires more sophisticated aggregation procedures than ERM.

We therefore consider the CAEW aggregate \eqref{equation:caew} and show that with a properly chosen temperature parameter $\beta$ it can achieve  faster rates than its ERM counterparts for $K < \sqrt{\frac{n}{\ln M}}$:
\begin{theorem} \label{theorem:caew-upper-bound}
    Consider the CAEW aggregate \eqref{equation:caew} over a dictionary $\cF$ with $\beta=c_1(B) K$ for $K < c_2(B) \sqrt{\frac{n}{\ln M}}$
    and $\beta=B\sqrt{\frac{2n}{\ln M}}$ otherwise. Then, for any $K,M,n$, and distribution $\pi$,
    \[
        \mathbb{E}_{\pi^{\otimes n}}[R(\hat{f}_{CAEW}] \leq \min_{1 \leq j \leq M}R(f_j) + c_3(B) \min\left\{\frac{K\ln(M)}{n}, \sqrt{\frac{\ln(M)}{n}}\right\},
    \] 
    where the constants $c_1(B), c_2(B), c_3(B)$ are given in the proof.

    \ignore{
        \[
            \mathbb{E}_{\pi^{\otimes n}}[R_{0/1}(\hat{f}_{CAEW}^{l,\cF})] - R_{0/1}^* \leq \sqrt{\min\{8e^2\frac{K\ln(M)}{n}, \sqrt{\frac{2\ln(M)}{n}}\}+2\bigl(\min_{1 \leq j \leq M}R(f_j)-R^* \bigr)}
        \] 
    }
\end{theorem}
In particular, in the binary case, CAEW achieves the optimal $\frac{\ln M}{n}$ convergence rate with respect to the logistic loss, matching the results for AEW \parencite{hogsgaard2026aggregationexponentialweightsoptimal}, $Q$-aggregation \parencite{lecuerigollet14Qaggregation} and the model selection aggregate proposed in \textcite{Rigollet_2012}. Moreover, as in the binary case, the CAEW outperforms the ERM counterparts with respect to multinomial logistic loss and is strictly superior  whenever $K < \sqrt{\frac{n}{\ln M}}$.

The following theorem shows that the rates established for CAEW in Theorem \ref{theorem:caew-upper-bound}, as well as those attained by the other aggregates mentioned above, cannot be improved, in the minimax sense, by any aggregate constrained to
$\operatorname{conv}(\cF)$. Therefore, the CAEW aggregate is rate-optimal across all proper aggregates within $\operatorname{conv}(\cF)$:
\begin{theorem} \label{theorem:proper-estimators-lower-bound}
    For every $K,M,n$, where $n>4\log_2(M)$ and $K \geq 3$, there exists a dictionary $\cF$ of size $M$ such that
    \[
    \inf_{\hat{f}\in \operatorname{conv}(\cF)}
    \sup_{\pi}
    \mathbb{E}_{\pi^{\otimes n}}\Big[R(\hat{f}) - \min_{1 \leq j \leq M} R(f_j)\Big] 
    \geq 
    c~ \min\left\{\frac{K\ln(M)}{n}, \sqrt{\frac{\ln(M)}{n}}\right\}
    \]
    for some $c>0$.
\end{theorem}

Importantly, the lower bound applies to all proper aggregates in $\operatorname{conv}(\cF)$, regardless of the aggregation procedure. 
Theorem \ref{theorem:proper-estimators-lower-bound} reveals an intriguing phase transition in aggregation rates when passing from the finite dictionary $\cF$ to its convex hull $\operatorname{conv}(\cF)$ in multiclass classification. 
Allowing aggregation over $\operatorname{conv}(\cF)$ 
indeed improves the rate when $K \lesssim \sqrt{\frac{n}{\ln(M)}}$.
Thus, the deterioration with $K$ is not merely a limitation of particular aggregation methods, but is unavoidable for proper aggregates in certain regimes. For larger $K$, the advantage of enlarging the aggregation class disappears, and the resulting rate $\sqrt{\frac{\ln(M)}{n}}$ coincides with that of Theorem \ref{theorem:erm-lower-bound}.

This limitation can be overcome by allowing improper aggregates outside $\operatorname{conv}(\cF)$. We now construct such an aggregate, which generally lies even outside $\operatorname{span}(\cF)$ and achieves the substantially faster rate $\frac{\ln(M)}{n}$, independently of $K$. \textcite{Rigollet_2012} established the optimality of this rate for model selection aggregation in binary classification ($K=2$) with fixed design. We further show in Theorem \ref{theorem:lower-bound-ms-general-case} that this rate remains optimal in multiclass classification for arbitrary $K$.

Let $p_j=\operatorname{softmax}(f_j),\;j=1,\ldots,M$, i.e.
$$
p_{jk}(x)=\frac{e^{f_{jk}(x)}}{\sum_{k'=1}^Ke^{f_{jk'}(x)}}
$$
which can be interpreted as class probabilities induced by $f_j$. Consider {\em probability-mixture aggregates} $f^P_\theta(x) \in  \mathbb{R}^K$, where
$$
    f^P_{\theta, k}(x) = \ln\left(\sum_{j=1}^M \theta_j p_{jk}(x)\right)=\ln\left(\sum_{j=1}^M \theta_j \frac{e^{f_{jk}(x)}}{\sum_{k'=1}^K e^{f_{jk'}(x)}}\right),\;\;\;k=1,\ldots,K
$$
and therefore, $\sum_{k=1}^K e^{f^P_{\theta,k}(x)}=1$.
The proposed procedure aggregates in the induced probability space  rather than in the original function space of $f_j$'s. Note that generally $f^P_\theta  \notin \operatorname{span}(\cF)$ and is an improper aggregate.

The corresponding multinomial logistic loss \eqref{equation:multinomial_loss}
is then
\begin{equation} \label{equation:multinomial_loss_probabilities}
    l(f^P_\theta(x),y)
    =
    -f^P_{\theta, y}(x)
    =
    -\ln\left(\sum_{j=1}^M \theta_j p_{jy}(x)\right)
    =
    -\ln\left(\sum_{j=1}^M \theta_j \frac{e^{f_{jy}(x)}}{\sum_{k=1}^K e^{f_{jk}(x)}}\right)
\end{equation}
Finally, consider the CAEW probability-mixture aggregate 
\begin{equation} \label{equation:probability_aggregate}
    \hat{f}^P_{CAEW}(x)=\ln\left(\sum_{j=1}^M \hat{\theta}_{CAEW,j} ~p_j(x)\right)
\end{equation}
with the CAEW weights $\hat{\theta}_{CAEW}$ from \eqref{equation:caew} for $\beta=1$, and $\ln( \cdot)$ is applied componentwise.

\begin{theorem}\label{theorem:optimal-ms-estimator}
Consider the CAEW probability-mixture aggregate \eqref{equation:probability_aggregate} over a dictionary $\cF$ with respect to multinomial logistic loss.
    For every $K,M,n$, and distribution $\pi$ 
    \[
        \mathbb{E}_{\pi^{\otimes n}}[R(\hat{f}^P_{CAEW})]  \leq \min_{1 \leq j \leq M}R(f_j) + \frac{\ln(M)}{n}
    \]
    \end{theorem}
The rate established in Theorem \ref{theorem:optimal-ms-estimator} does not depend on the number of classes. We now show that this is the optimal rate for the multinomial logistic loss.

\begin{theorem}\label{theorem:lower-bound-ms-general-case}
    For every $K,M$ and $n \geq \log_2(M)$, there exists a dictionary $\cF$ of size $M$ such that 
    \begin{gather}
        \inf_{\hat{f}}
        \sup_{\pi}
        \mathbb{E}_{\pi^{\otimes n}}\Big[R(\hat{f}) - \min_{1 \leq j \leq M} R(f_j)\Big] 
        \geq 
        c\frac{\ln(M)}{n}
    \end{gather}
    for some $c>0$.
\end{theorem}

\begin{corollary}\label{corollary:optimal-estimator-both-losses}
    The CAEW probability-mixture aggregate $\hat{f}^P_{CAEW}$ in Theorem \ref{theorem:optimal-ms-estimator} is rate-optimal with respect to the multinomial logistic loss, with optimal rate $\frac{\ln(M)}{n}$ independent of $K$.
\end{corollary}
From Theorems \ref{theorem:proper-estimators-lower-bound}-\ref{theorem:lower-bound-ms-general-case} it follows that unless the size of the dictionary grows exponentially with $n$, the original CAEW aggregate \eqref{equation:caew} defined on $\cF$ is rate-optimal with respect to the multinomial logistic loss only for very small $K=\mathcal O(1)$. In all other (``truly multiclass'') regimes, attaining optimal rates requires moving beyond $\operatorname{conv}(\cF)$. The key reason lies in the geometry of the multinomial logistic loss. Its exp-concavity decreases at the rate $\frac{1}{K}$ with the number of classes $K$ when parameterized in score space by the functions $f_j$, whereas in the induced probability space the logistic loss is 1-exp-concave (see the proof of Theorem \ref{theorem:optimal-ms-estimator}). Consequently, the proposed probability-mixture aggregate bypasses this loss of curvature and achieves the optimal aggregation rate independently of $K$.

\ignore{
\newcolumntype{Y}{>{\centering\arraybackslash}X}

The results of this section are summarized in Table \ref{table:results-ms}.

\begin{table}[htbp]
    \centering
    \small
    \begin{tabularx}{\textwidth}{|Y|Y|Y|}
        \hline
        \begin{gather*}
            \hat{f}
        \end{gather*}
        & \begin{gather*}
            \mathbb{E}_{\pi^{\otimes n}}[R(\hat{f})] - \min_{1 \leq j \leq M}R(f_j)
        \end{gather*}
        & \begin{gather*}
            \mathbb{E}_{\pi^{\otimes n}}[R_{0/1}(\hat{f})] - R^*
        \end{gather*} \\
        \hline
        \begin{gather*}
            \hat{f}_{ERM}^{l,\cF}
        \end{gather*}
        & \begin{gather*}
            \sqrt{\frac{\ln(M)}{n}}
        \end{gather*}
        & \begin{gather*}
           \Bigg(
            32\sqrt{\frac{\ln(M)}{n}}
            \\+ 2\bigl(\min_{1 \leq j \leq M}R(f_j)-R^* \bigr)
            \Bigg)^{1/2} 
        \end{gather*} \\
        \hline

        $\hat{f}_{ERM}^{l_{0/1},\cF}$ 
        & -- 
        & \begin{gather*}
            \sqrt{\frac{\ln(M)}{n}} \\
            \text{compared to} \;\; \min_{1 \leq j \leq M}R_{0/1}(f_j)
        \end{gather*} \\
        \hline
        \begin{gather*}
            \hat{f}_{CAEW}^{l,\cF}
        \end{gather*}
        & \begin{gather*}
            \min\left\{
            \frac{K\ln(M)}{n},
            \sqrt{\frac{\ln(M)}{n}}
            \right\}
        \end{gather*}
        & \begin{gather*}
            \Bigg(
            \min\left\{
            8e^2\frac{K\ln(M)}{n},
            \sqrt{\frac{2\ln(M)}{n}}
            \right\}
            \\+ 2\bigl(\min_{1 \leq j \leq M}R(f_j)-R^* \bigr)
            \Bigg)^{1/2}
        \end{gather*} \\
        \hline
        \begin{gather*}
            \hat{f}_{Prob\text{-}CAEW}^{l,\cF}
        \end{gather*}
        (Theorem \ref{theorem:optimal-ms-estimator}) 
        & \begin{gather*}
            \frac{\ln(M)}{n}
        \end{gather*}
        & \begin{gather*}
            \Bigg(
            2\frac{\ln(M)}{n} \\
            + 2\bigl(\min_{1 \leq j \leq M}R(f_j)-R^* \bigr)
            \Bigg)^{1/2}
        \end{gather*} \\
        \hline
    \end{tabularx}
    \caption{MS-type aggregation convergence rates summarized}
    \label{table:results-ms}
\end{table}
}

\section{Convex aggregation} \label{section:convex}
In this section, we consider convex aggregation, where $\Theta=\Delta_M$ and $\cF_\Theta=\operatorname{conv}(\cF)$. For a given ($0/1$ or multinomial logistic) loss the ultimate goal is to construct an aggregated classifier $\hat{f}$ based on the training sample $\mathcal S^n$ satisfying
$$    
  \mathbb{E}_{\pi^{\otimes n}}\big[R(\hat f)\big] 
    \leq 
    \inf_{\theta \in \Delta_M} R(f_\theta) 
    +
    \epsilon(M,n,\Pi)
$$
with the minimal (rate-optimal) remainder term $\epsilon(M,n,\Pi)$. The optimal rates for convex aggregation cannot, in general, be faster than those for model selection aggregation. 

We start by considering the ERM aggregate over $\operatorname{conv}(\cF)$ in \eqref{equation:erm-convex-hull} with respect to the $0/1$ loss.

\begin{theorem}\label{theorem:convex-erm-0/1}
    There exists $c>0$ such that for every $K,M,n,$ dictionary $\cF$, and distribution $\pi$, the ERM aggregate $\hat{f}_{ERM}^{l_{0/1},\operatorname{conv}(\cF)}$ from \eqref{equation:erm-convex-hull} satisfies:
     \[
     \mathbb{E}_{\pi^{\otimes n}}[R_{0/1}(\hat{f}_{ERM}^{l_{0/1},\operatorname{conv}(\cF)})] 
     \leq  
     \inf_{f \in \operatorname{conv}(\cF)}R_{0/1}(f) 
     + 
     c\min\left\{1,\sqrt{\frac{M \ln(K)}{n}}\right\}
     \]
\end{theorem}
From \textcite{daniely2014optimallearnersmulticlassproblems} it follows that the $\ln(K)$-factor in the upper bound for the ERM aggregate is generally unavoidable.
On the other hand, from the VC theory for convex linear classifiers, it is known that for binary classification
\[
\sup_{\cF: |\cF| = M}
\inf_{\hat f\in\operatorname{conv}(\cF)}
\sup_{\pi}
\mathbb{E}_{\pi^{\otimes n}}
\left[
R_{0/1,\pi}(\hat f)
-
\inf_{f\in\operatorname{conv}(\cF)}
R_{0/1,\pi}(f)
\right]
\gtrsim
\min \Big\{1, \sqrt{\frac{M}{n}}\Big\}.
\]

Since multiclass classification contains binary classification as a particular case, the same lower bound applies to the multiclass setting. Therefore, Theorem~{\ref{theorem:convex-erm-0/1}} implies that the ERM aggregate
\(\hat f_{\mathrm{ERM}}^{\,l_{0/1},\operatorname{conv}(\cF)}\)
is nearly rate-optimal over \(\operatorname{conv}(\cF)\), up to an additional
\(\sqrt{\ln K}\)-factor.
Unfortunately, optimization of the empirical $0/1$ risk over $\operatorname{conv}(\cF)$ is nonconvex and generally computationally infeasible, although the estimator remains of theoretical interest. 

We now consider the ERM aggregate \eqref{equation:erm-convex-hull} with respect to the multinomial logistic loss. Since the loss is convex, the corresponding ERM can be computed efficiently.

\begin{theorem}\label{theorem:erm-convex-hull-upper-bound}
    For every $K,M,n$, dictionary $\cF$, and distribution $\pi$, the aggregate $\hat{f}_{ERM}^{l,\operatorname{conv}(\cF)}$ from \eqref{equation:erm-convex-hull} satisfies:
     \[
         \mathbb{E}_{\pi^{\otimes n}}[R(\hat{f}_{ERM}^{l,\operatorname{conv}(\cF)})] \leq \inf_{f \in \operatorname{conv}(\cF)}R(f) + c(B) \min\left\{\sqrt{\frac{\ln(M)}{n}}, \frac{KM}{n}\right\}
     \]
        where $c(B)>0$ depends only on $B$.
\end{theorem}

The following theorem establishes a lower bound for convex aggregation:
\begin{theorem}\label{theorem:convex-aggregation-lower-bound}
    Assume that $\mathcal{X}$ contains at least $2^M$ distinct points.
    For every $M$ and $n$  such that $\ln(1+\frac{M}{\sqrt{n}})<n$, and every $K$ there exists a dictionary $\cF$ of size $M$ such that 
    \begin{gather} \label{equation:lower_bound_convex}
        \inf_{\hat{f}} 
        \sup_{\pi}
        \mathbb{E}_{\pi^{\otimes n}}\Big[R_{\pi}(\hat{f}) - \inf_{f \in \operatorname{conv}(\cF)}R_{\pi}(f) \Big] 
        \geq 
        c\min\left\{\sqrt{\frac{\ln(1+\frac{M}{\sqrt{n}})}{n}},\frac{M}{n}\right\}
\end{gather}
for some $c>0$.
\end{theorem}
\noindent
As in model selection aggregation, the lower bound \eqref{equation:lower_bound_convex} does not depend on the number of classes and coincides with that obtained in \textcite{Rigollet_2012} for binary classification with fixed design.

Theorems \ref{theorem:erm-convex-hull-upper-bound} and \ref{theorem:convex-aggregation-lower-bound} imply that the ERM aggregate $\hat{f}_{ERM}^{l,\operatorname{conv}(\cF)}$ from \eqref{equation:erm-convex-hull} is nearly rate-optimal (up to differences in the log-term) over $\operatorname{conv}(\cF)$ when either $K = \Omega(\frac{\sqrt{n \ln M}}{M})$ and $M = \Omega(\sqrt{n})$, or $K=\mathcal{O}(1)$ and $M = \mathcal{O}(\sqrt{n})$.  
In the intermediate regime, however, a substantial gap involving a linear dependence on $K$ remains between the upper and lower bounds.

In fact, the linear dependence on $K$ in some regimes in Theorem \ref{theorem:erm-convex-hull-upper-bound} is inherent not only to ERM but to {\em any proper} convex aggregate. 
Indeed, since $\inf_{f \in \operatorname{conv}\cF}R(f) \leq \min_{f \in \cF}R(f)$, 
the lower bound of Theorem~\ref{theorem:proper-estimators-lower-bound}, established for all aggregates in $\operatorname{conv}(\mathcal F)$ relative to the model-selection oracle, applies {\em a fortiori} when the excess risk is measured relative to the convex-aggregation oracle.

Hence, similar to model selection, to achieve optimal rates for convex aggregation one has to apply improper aggregates. In particular, consider the {\em continuous} version of CAEW probability-mixture $\hat{f}^P_{CAEW}$ with $\beta=1$ defined in \eqref{equation:probability_aggregate}.
Namely,
\begin{equation} \label{equation:probability_aggregate_continuous}
    \hat{f}^{cP}_{CAEW}(x)=\ln \int_{\Delta_M} \operatorname{softmax}(f_\theta(x)) \widehat{w}_{CAEW}(\theta)d\theta,
\end{equation}    
where
$$
f_\theta(x)=\sum_{j=1}^M \theta_j f_j(x),\;\;\;\;
\widehat{w}_{CAEW}(\theta)=\frac{1}{n}\sum_{i=0}^{n-1} \frac{e^{-i L^i (f_\theta)}}{\int_{\Delta_M}e^{-i L^i (f_\theta)}d\theta},\;\;\;\;L^0(f_\theta)=0,
$$
\(d\theta\) denotes the \((M-1)\)-dimensional volume measure on \(\Delta_M\) 
and $\ln( \cdot)$ is applied componentwise. 

Unlike the original $\hat{f}^P_{CAEW}$ in \eqref{equation:probability_aggregate} used for model selection aggregation, the proposed aggregate in \eqref{equation:probability_aggregate_continuous} averages the induced probabilities over the entire simplex $\Delta_M$, rather than over its vertices. This extension is motivated by the fact that the oracle risk for convex aggregation is defined as the minimum risk over $\operatorname{conv}(\cF)$, rather than over the finite dictionary $\cF$.

We have the  following oracle inequality for $\hat{f}^{cP}_{CAEW}$:
\begin{theorem}\label{theorem:continuous-CAEW-logistic-upper-bound}
    For every $K, n, M \geq 2$, dictionary $\cF$, and distribution $\pi$, 
     \[
         \mathbb{E}_{\pi^{\otimes n}}[R(\hat{f}_{CAEW}^{cP})] \leq 
         \inf_{f \in \operatorname{conv}(\cF)}R(f) + \frac{M-1}{n}\Bigl(1+\ln(1+\frac{4B n}{M-1})\Bigr)
     \]    
\end{theorem}

Theorems \ref{theorem:erm-convex-hull-upper-bound} and \ref{theorem:continuous-CAEW-logistic-upper-bound} 
suggest an optimal aggregation strategy for convex aggregation under the multinomial logistic loss  that switches between ERM and continuous probability-mixture CAEW, depending on the relative values of $M$ and $n$. This yields the following corollary:
\begin{corollary} \label{corollary:logistic-hull-upper-bound} Let $\hat{f}^{l,\operatorname{conv}(\cF)}=\hat{f}^{l,\operatorname{conv}(\cF)}_{ERM}$ from \eqref{equation:erm-convex-hull} if $\min\left\{\sqrt{\frac{\ln(M)}{n}},\frac{K}{M}\right\} \leq \frac{M}{n}\max\left\{1,\ln\Bigl(\frac{n}{M}\Bigr)\right\}$
and $\hat{f}^{l,\operatorname{conv}(\cF)}=\hat{f}_{CAEW}^{cP}$ from \eqref{equation:probability_aggregate_continuous} otherwise. Then,
\[
\mathbb E_{\pi^{\otimes n}}
\left[R\left(\hat f^{l,\operatorname{conv}(\cF)}\right)\right]
\le
\inf_{f\in\operatorname{conv}(\cF)}R(f)
+c(B)\min\left\{
\sqrt{\frac{\ln(M)}{n}},
\frac{KM}{n},
\frac{M}{n}\max\left\{1,\ln\left(\frac nM\right)\right\}
\right\}
\]
\ignore{ \[
         \mathbb{E}_{\pi^{\otimes n}}[R(\hat{f}^{l,\operatorname{conv}(\cF)})] \leq 
         \inf_{f \in \operatorname{conv}(\cF)}R(f) + c(B) \min\left\{\sqrt{\frac{\ln(M)}{n}},\frac{M}{n}\max\{1,\ln\Bigl(\frac{n}{M}\Bigr)\}\right\}
     \] 
 }    
and, therefore, it achieves the optimal rate \eqref{equation:lower_bound_convex} (independent of $K$), up to an extra log-factor in certain regimes.
\end{corollary}
Thus, as in model selection aggregation, allowing improper aggregates by averaging the induced class probabilities and mapping the result back to score-function space yields nearly optimal convergence rates independently of $K$. The additional $\frac{KM}{n}$ term provides a sharper upper bound when $K$ is sufficiently small relative to $M$ and $n$, in which case proper ERM aggregation is preferable.

\ignore{
In Section~\ref{section:ms-aggregation}, we showed that although the optimal rate for model selection aggregation with respect to the logistic loss does not depend on the number of classes, the lower bound on the convergence rate of proper aggregates within $\operatorname{conv}(\cF)$ does.
This observation suggests that the dependence on $K$ may be an inherent limitation of proper aggregation, leading to the conjecture that, as in model selection, it can be avoided in the convex case by using improper aggregates.
}

\ignore{
The results of this section are summarized in table \ref{table:results-convex}.
\begin{table}[htbp]
    \centering
    \small
    \begin{tabularx}{\textwidth}{|Y|Y|Y|}
        \hline
        \begin{gather*}
            \hat{f}
        \end{gather*}
        & \begin{gather*}
            \mathbb{E}_{\pi^{\otimes n}}[R(\hat{f})] - \inf_{f \in \operatorname{conv}(\cF)}R(f)
        \end{gather*}
        & \begin{gather*}
            \mathbb{E}_{\pi^{\otimes n}}[R_{0/1}(\hat{f})] - R^{*}
        \end{gather*} \\
        \hline
        \begin{gather*}
            \hat{f}_{ERM}^{l,\operatorname{conv}(\cF)}
        \end{gather*}
        & \begin{gather*}
            \min\{\sqrt{\frac{\ln(M)}{n}}, \frac{KM}{n}\}
        \end{gather*}
        & \begin{gather*}
           (c * \min\{\sqrt{\frac{\ln(M)}{n}}, \frac{KM}{n}\} \\
           + \bigl(\inf_{f \in \operatorname{conv}(\cF)}R(f)-R^*\bigr))^{\frac{1}{2}}
        \end{gather*} \\
        \hline
        \begin{gather*}
            \hat{f}_{ERM}^{l_{0/1},\operatorname{conv}(\cF)}
        \end{gather*}
        & \begin{gather*}
            -
        \end{gather*} 
        & \begin{gather*}
            \sqrt{\frac{M\ln(K)}{n}} \\ \text{compared to} \;\; \inf_{f \in \operatorname{conv}(\cF)}R_{0/1}(f)
        \end{gather*} \\
        \hline
        \begin{gather*}
            \hat{f}_{Prob-CAEW}^{l,\operatorname{conv}(\cF)} \;\; (M>\sqrt{n})
        \end{gather*}
        (Theorem \ref{theorem:optimal-convex-stronger-assumption}) 
        & \begin{gather*}
            \sqrt{\frac{\ln(1+\frac{M}{\sqrt{n}})}{n}}
        \end{gather*}
        & \begin{gather*}
            \Bigg(c\sqrt{\frac{\ln(1+\frac{M}{\sqrt{n}})}{n}} \\
            + 2\bigl(\min_{f \in \operatorname{conv}(\cF)}R(f)-R^*\bigr)\Bigg)^{\frac{1}{2}}
        \end{gather*} \\
        \hline
    \end{tabularx}
    \caption{MS-type aggregation convergence rates summarized}
    \label{table:results-convex}
\end{table}
}

\section{Aggregation in probability space} \label{section:probaggregation}
In Sections \ref{section:ms-aggregation} and \ref{section:convex}, we constructed improper aggregates by averaging the induced class probabilities and then mapping the result back to score-function space, with the oracle defined in the latter. 
We now consider a different aggregation framework in which both the aggregation class and the oracle are defined directly in probability space.

\ignore{
As discussed in the previous sections, each score function $f_j \in \cF$ induces a probability vector $p_j(x)=\operatorname{softmax}(f_j(x)) \in \mathbb{R}^K$. Since classification ultimately depends on these induced class probabilities, it is natural to investigate aggregation of the probability vectors $p_j$ rather than the original score functions $f_j$. 
}

Consider now the multinomial logistic loss in terms of the induced probabilities:  
$$
l_{prob}(p(x),y)=-\ln \Big(p_y(x) \Big)
$$
(see \eqref{equation:multinomial_loss_probabilities})
and the corresponding risk $R_{prob}(p)=\mathbb{E}_\pi\Big[ l_{prob}(p(X),Y)\Big]$.

Let $\cP = \{p_1, \dots p_M\}$, where $p_j=\operatorname{softmax}(f_j),\; j=1,\ldots,M$. Recall that we assume that $\lVert f_j(x) \rVert_{\infty} \leq B$, and therefore $\frac{e^{-2B}}{K} \leq p_{jk}(x) \leq \frac{e^{2B}}{K},$ for all $\;j=1,\ldots, M, k=1, \dots, K$, and $x \in {\cal X}$. 

For a given $\Theta \subseteq \Delta_M$ consider the class
$\mathcal{P}_{\Theta} = \{p_{\theta}=\sum_{j=1}^{M}\theta_j p_j,\;\theta\in\Theta\}$.
The goal is to construct an aggregate $\hat{p}$ whose risk is as close as possible to $\inf_{p \in \cP_{\Theta}} R_{prob}(p)$, i.e.,
satisfying the oracle inequality
\[
    \mathbb{E}_{\pi^{\otimes n}}\big[R_{prob}(\hat p)\big] 
    \leq 
    \inf_{p \in \cP_\Theta} R_{prob}(p) 
    +
    \epsilon(M,n,\Pi)
\]
with the smallest possible remaining term.

\ignore{
We already introduced the idea of aggregation in probability space  in Sections  \ref{section:ms-aggregation} and \ref{section:convex} in the context of model selection and convex aggregation, where we proposed the probability-mixture CAEW aggregate \eqref{equation:probability_aggregate}. This aggregate is obtained by applying the CAEW procedure defined in~\eqref{equation:caew} to the induced class probabilities and then mapping the resulting probability mixture back to the score-function space. Now we are interested in  aggregation of probabilities themselves.
}

We start from model selection aggregation, where
$\Theta = \{e_1, \dots e_M \}$. In this case
\begin{gather}\label{equation:equivalence-loss-on-vertices}
    \min_{1 \leq j \leq M}R_{prob}(p_j)
    =
    \min_{1 \leq j \leq M} \mathbb{E}_\pi\Big[ -\ln \Big(p_{jY}(X) \Big)\Big] 
    =
    \min_{1 \leq j \leq M} \mathbb{E}_\pi\Big[ -\ln \Big(\frac{e^{f_{jY}(X)}}{\sum_{k=1}^Ke^{f_{jk}(X)}}\Big)\Big] 
    =
    \min_{1 \leq j \leq M}R(f_j)
\end{gather}
and, therefore, model selection probability-space aggregation is equivalent to the model selection aggregation of score functions $f_j$ considered in Section \ref{section:ms-aggregation}. In addition, the CAEW aggregate $\hat{p}_{CAEW}^{l,\cP}$ with respect to the loss $l_{prob}$ defined in \eqref{equation:caew} is equivalent to the CAEW probability-mixture aggregate \eqref{equation:probability_aggregate} in the following sense:

\begin{equation}\label{equation:equivalence-between-losses}
    l_{prob}(\hat{p}_{CAEW}^{l,\cP}(x),y)
    =
    -\ln\Big(\hat{p}_{CAEW,y}^{l,\cP}(x)\Big) 
    = 
    -\hat{f}_{CAEW,y}^{P}(x)
    =
    l\Big(\hat{f}_{CAEW}^{P}(x), y\Big),
\end{equation} 
where the last equality follows from \eqref{equation:multinomial_loss_probabilities}. Equations~\eqref{equation:equivalence-loss-on-vertices} and~\eqref{equation:equivalence-between-losses}, together with Theorem~\ref{theorem:optimal-ms-estimator}, imply the following corollary:
\begin{corollary}\label{corollary:optimal-rate-probability-mixture-ms-aggregation}
    The optimal rate for probability-space model selection aggregation with respect to the logistic loss is $\frac{\ln M}{n}$ and it is achieved by $\hat{p}_{CAEW}^{l,\cP}$ with a temperature $\beta=1$.
\end{corollary}

We now consider convex probability-space aggregation, where $\Theta=\Delta_M$ or equivalently, $\cP_\Theta=\operatorname{conv}(\cP)$. Generally, in this case, $\inf_{\theta \in \Theta}R(f_{\theta}) \neq \inf_{\theta \in \Theta}R_{prob}(p_{\theta})$ and the optimal rates for score-function aggregation and probability-space aggregation may differ.

Consider the ERM aggregate over $\operatorname{conv}(\cP)$ with respect to $l_{prob}$:  
\begin{equation} \label{equation:p-aggregation}
    \hat{p}_{ERM}^{l_{prob},\operatorname{conv}(\cP)} 
    = 
    \operatorname*{arg\,min}_{p \in \operatorname{conv}(\cP)} \frac{1}{n}\sum_{i=1}^n l_{prob}(p(x_i),y_i)
\end{equation}
We first derive the upper bound for its risk:
\begin{theorem}\label{theorem:probability-mixture-convex-aggregation-erm}
    For every $K,M,n$, dictionary $\cF$, and distribution $\pi$, the aggregate $\hat{p}_{ERM}^{l_{prob},\operatorname{conv}(\cP)}$ from \eqref{equation:p-aggregation} satisfies
     \[
         \mathbb{E}_{\pi^{\otimes n}}[R_{prob}(\hat{p}_{ERM}^{l_{prob},\operatorname{conv}(\cP)})] 
         \leq 
         \inf_{p \in \operatorname{conv}(\cP)}R_{prob}(p) 
         + 
         c(B) \min\left\{\sqrt{\frac{\ln(M)}{n}}, \frac{M}{n}\right\}
     \]
     where $c(B)>0$ depends only on $B$.
\end{theorem}

Hence, unlike the ERM aggregate $\hat{f}_{ERM}^{l,\operatorname{conv}(\cF)}$ of score functions from Section \ref{section:convex}, the bound for the probability-space ERM $\hat{p}_{ERM}^{l_{prob},\operatorname{conv}(\cP)}$ established in Theorem \ref{theorem:probability-mixture-convex-aggregation-erm} does not depend on $K$. We now show that it is essentially tight (up to a logarithmic gap), which establishes the near-optimality of the $\hat{p}_{ERM}^{l_{prob},\operatorname{conv}(\cP)}$ aggregate.

\begin{theorem}\label{theorem:probability-mixture-convex-aggregation-lower-bound} 
     Assume that $\mathcal{X}$ contains at least $2^M$ distinct points. 
    For every $M$ and $n$  such that $\ln(1+\frac{M}{\sqrt{n}})<n$, and every $K$ 
    there exists a dictionary $\cF$ of size $M$ such that
    \begin{gather}
        \inf_{\hat{p}}\sup_\pi \mathbb{E}_{\pi^{\otimes n}}\Big[R_{prob}(\hat{p}) - \inf_{p \in \operatorname{conv}(\cP)}R_{prob}(p)\Big] 
        \geq 
        c\min\left\{\sqrt{\frac{\ln(1+\frac{M}{\sqrt{n}})}{n}},\frac{M}{n}\right\}
\end{gather}
for some $c>0$.
\end{theorem}

\ignore{
The results of this section show that aggregating the induced probabilities rather than the original score functions removes the dependence on the number of classes and achieves a near-optimal aggregation rate (up to log factors). As already discussed in Section \ref{section:ms-aggregation}, the key reason lies in the geometry of the multinomial logistic loss.
Under the usual score-space parameterization, the exp-concavity constant deteriorates as $\frac{1}{K}$ with the number of classes, whereas in the induced probability space the logistic loss is 1-exp-concave. Consequently, probability-mixture aggregation avoids this deterioration in curvature and attains the near-optimal aggregation rate regardless of $K$.
}
Theorems \ref{theorem:probability-mixture-convex-aggregation-erm} and \ref{theorem:probability-mixture-convex-aggregation-lower-bound} establish that ERM over the convex hull of the induced probability vectors achieves the minimax-optimal convergence rate, up to logarithmic factors, independently of the number of classes $K$. Unlike convex aggregation in the original score-function space, this result is achieved by a single, computationally tractable ERM procedure, without requiring improper aggregation or switching between different procedures.

\vspace{.2cm}
\noindent
{\bf Acknowledgments}. 
The work was supported by the Israel Science Foundation (ISF), Grant ISF 1095/22.

\section{Appendix}

\begin{proof}[\textbf{Proof of Lemma~\ref{lemma:loss-strongly-convex-and-smooth-wrt-function}}]
For a given $f(x) \in \mathbb{R}^K$ let $p(x)=\operatorname{softmax}(f(x))$ denote the corresponding induced multinomial probabilities.
A straightforward calculus shows that the Hessian of the multinomial logistic loss is a symmetric positive semidefinite matrix $V(f(x))=diag(p(x))-p(x) p(x)^T$. Therefore, since the multinomial logistic loss is twice continuously differentiable, by Taylor's expansion
$$
l(f_1(x),y)-l(f_2(x),y) = \nabla l(f_2(x),y)^T(f_1(x)-f_2(x)) +\frac{1}{2} (f_1(x)-f_2(x))^T V(\tilde{f}(x))(f_1(x)-f_2(x)),
$$
where $\tilde{f}(x)=\alpha f_2(x)+(1-\alpha)f_{1}(x)$ for some $0 \leq \alpha \leq 1$.

Let $\lambda_{(1)}(x)  \leq \ldots \leq \lambda_{(K)}(x)$ be the eigenvalues of $V(\tilde{f}(x))$.
Clearly, $\lambda_{(1)}=0$ and the null space $\operatorname{Null}(V(\tilde{f}))=\operatorname{span}\left((1,\ldots,1)^T\right)$. By the symmetric constraint \eqref{equation:symmetric_constraint}, $f_1(x)-f_2(x) \in \operatorname{Null}(V(\tilde{f}))^\perp$. 
Therefore,
$$
\lambda_{(2)} \lVert f_1(x)-f_2(x) \rVert^2_2 \leq  (f_1(x)-f_2(x))^T V(\tilde{f}(x))(f_1(x)-f_2(x)) \leq \lambda_{(K)} \lVert f_1(x)-f_2(x) \rVert^2_2
$$
which implies
\begin{equation} \nonumber
\begin{split}
\nabla & l(f_2(x),y)^T(f_1(x)-f_2(x)) + \frac{\lambda_{(2)}}{2} \lVert f_1(x)-f_2(x) \rVert^2_2 \\
& \leq l(f_1(x),y)-l(f_2(x),y) \\
& \leq \nabla l(f_2(x),y)^T(f_1(x)-f_2(x)) + \frac{\lambda_{(K)}}{2} \lVert f_1(x)-f_2(x) \rVert^2_2
\end{split}
\end{equation}
On the other hand, it is known that 
$\tilde{p}_{(1)}(x) \leq  \lambda_{(2)}(x) \leq \ldots \leq \lambda_{(K)}(x) \leq \tilde{p}_{(K)}(x)$ \parencite{10.1111/j.2517-6161.1996.tb02081.x}. 

Finally, since  $\lVert f_1(x) \rVert_\infty \leq B$ and $\lVert f_2(x) \rVert_\infty \leq B$, we have $\lVert \tilde{f}(x) \rVert_\infty \leq B$ and therefore $\frac{e^{-2B}}{K} \leq \tilde{p}_{(1)}(x) \leq \tilde{p}_{(K)}(x) \leq \frac{e^{2B}}{K}$. 

Combining the above inequalities completes the proof.
\end{proof}

\begin{proof}[\textbf{Proof of Lemma~\ref{lemma:logistic-loss-geometry}}]
    Consider $\theta_1, \theta_2 \in \Delta_M$. Similar to the proof of Lemma~\ref{lemma:loss-strongly-convex-and-smooth-wrt-function}, by Taylor's expansion 
    \begin{equation*}
    \begin{split}
        l(F(x)\theta_1,y)& -l(F(x)\theta_2,y)  
        \\
        & =\nabla_\theta l(F(x)\theta_2,y)^T(\theta_1-\theta_2)
        + 
        \frac{1}{2}(\theta_1-\theta_2)^TF(x)^T \left(\operatorname{diag}(p_{\tilde \theta}(x))-p_{\tilde \theta}(x) p_{\tilde \theta}(x)^T\right)F(x)(\theta_1-\theta_2),
    \end{split}
    \end{equation*}
        where $\tilde{\theta}=\alpha \theta_2+(1-\alpha) \theta_1,\;0 \leq \alpha \leq 1$ and $p_{\tilde \theta}(x)=\operatorname{softmax}(F(x)\tilde{\theta})$. 
    
    By the symmetric constraint \eqref{equation:symmetric_constraint}, $F(x)(\theta_1-\theta_2) \in \operatorname{span}((1,\ldots,1)^T)^\perp$, and one again has
    \begin{equation} \nonumber
    \begin{split}
    \nabla_\theta & l(F(x)\theta_2,y)^T(\theta_1-\theta_2) + \frac{e^{-2B}}{2K}(\theta_1-\theta_2)^TF(x)^TF(x)(\theta_1-\theta_2) \\
    & \leq l(F(x)\theta_1,y) - l(F(x)\theta_2,y) \\
    & \leq  \nabla_\theta l(F(x)\theta_2,y)^T(\theta_1-\theta_2) + \frac{e^{2B}}{2K}(\theta_1-\theta_2)^TF(x)^TF(x)(\theta_1-\theta_2)
    \end{split}
    \end{equation}
    Thus, the loss $l(f_\theta(x),y)$ is $\frac{e^{-2B}}{K}$-strongly convex and $\frac{e^{2B}}{K}$-smooth with respect to the seminorm induced by $H(x)=F(x)^TF(x)$. 
    
    Furthermore, note that $\nabla_\theta l(F(x)\theta,y)=F(x)^T(p_\theta(x)-e^y)$,
    where $p_\theta(x)=\operatorname{softmax}(F(x)\theta)$. Hence, 
    \begin{gather*}
        |l(F(x)\theta_1,y)-l(F(x)\theta_2,y)| =|\nabla_\theta l(F(x){\tilde \theta},y)^T(\theta_1-\theta_2)| 
        = \\ 
        |(p_\theta(x)-e_y)^TF(x)(\theta_1-\theta_2)| \leq 
        \lVert p_{\tilde \theta}(x)-e_y \rVert_2~ \lVert \theta_1-\theta_2 \rVert_H 
    \end{gather*}
    for some $\tilde{\theta}=\alpha \theta_1+(1-\alpha)\theta_2,\; 0 \leq \alpha \leq 1$, where
    $$
    \lVert p_{\tilde \theta}(x)-e_y \rVert^2_2=(1-p_{\tilde \theta, y}(x))^2+\sum_{k \neq y} (p_{\tilde \theta, k}(x))^2 \leq (1-p_{\tilde \theta, y}(x))^2+\Big(\sum_{k \neq y} p_{\tilde \theta, k}(x)\Big)^2 = 2(1-p_{\tilde \theta, y}(x))^2 \leq 2
    $$
    Thus, the loss $l(F(x)\theta,y)$ is $\sqrt{2}$-Lipschitz with respect to the seminorm induced by $H$. 
    
    Lastly, Proposition 3.2 in \textcite{puchkin2023exploringlocalnormsexpconcave} implies that $\sqrt{2}$-Lipschitz and $\frac{e^{-2B}}{K}$-strong-convex function is $\frac{1}{2e^{2B}K}$-exp-concave.

\end{proof}

\begin{lemma}\label{lemma:discrete-kl-lower-bound}
For any two probability vectors $p(x),q(x)\in\Delta_K$, any
$A\subseteq[K]$, and any distribution of $X$,
\[
\mathbb{E}_X\!\left[
\operatorname{KL}(p(X) \parallel q(X))
\right]
\ge
2\mathbb{E}_X\!\left[
\bigl(p_A(X)-q_A(X)\bigr)^2
\right],
\]
where
\[
p_A(x)=\sum_{k\in A}p_k(x),
\qquad
q_A(x)=\sum_{k\in A}q_k(x).
\]
   
\end{lemma}

\begin{proof}[\textbf{Proof of Lemma~\ref{lemma:discrete-kl-lower-bound}}]
    For a given $x \in \mathcal{X}$ we have
    \begin{gather*}
        KL(p(x) \parallel q(x)) = \sum_{k=1}^K p_k(x)\ln\frac{p_k(x)}{q_k(x)} 
        = 
        \sum_{k \in A} p_k(x)\ln\frac{p_k(x)}{q_k(x)}  
        + 
        \sum_{k \notin A} p_k(x)\ln\frac{p_k(x)}{q_k(x)}  
    \end{gather*}

    Using the log-sum inequality:
    \begin{gather*}
        \sum_{k \in A} p_k(x)\ln\frac{p_k(x)}{q_k(x)}  
        + 
        \sum_{k \notin A} p_k(x)\ln\frac{p_k(x)}{q_k(x)} 
        \geq
        p_A(x)\ln \frac{p_A(x)}{q_A(x)} + (1-p_A(x))\ln\frac{1-p_A(x)}{1-q_A(x)} \\
        = 
        KL\Big(Ber(p_A(x))\parallel Ber(q_A(x))\Big)
    \end{gather*}

    By Pinsker's inequality,
   $$
        KL\Big(Ber(p_A(x)) \parallel Ber(q_A(x))\Big) \geq 2 \lVert (Ber(p_A(x))-Ber(q_A(x)) \rVert^2_{TV}=(p_A(x)-q_A(x))^2
    $$
    Taking expectations with respect to $x$ completes the proof. 
\end{proof}

\begin{proof}[\textbf{Proof of Theorem~\ref{theorem:erm-upper-bound-expectation}}]
The proof of the theorem immediately follows from the general Proposition \ref{proposition:upper_bound} below by noting that $l_{0/1}$-loss is obviously bounded on the interval $[0,1]$, while for the multinomial logistic loss, since $\lVert f_j(x) \rVert _\infty\le B$,
\begin{equation} \label{equation:bounded_loss}
-2B+\ln K
\le
l(f_j(x),y)
\le
2B+\ln K
\end{equation}
implying that the loss is bounded on an interval of length $4B$.
\end{proof}

\begin{proposition} \label{proposition:upper_bound} Consider a general loss $l: \mathbb{R}^K \times {\cal Y} \to \mathbb{R}$ bounded on an interval of length $h$. Then,
$$
\mathbb{E}_{\pi^{\otimes n}}
\Bigl[
R\bigl(\hat f_{\mathrm{ERM}}^{l,\cF}\bigr)\Bigr]
-\min_{1 \leq j \leq M}R(f_j)
\leq
4h\sqrt{\frac{\ln(M)}{2n}}
$$
\end{proposition}
\begin{proof}[\textbf{Proof of Proposition~\ref{proposition:upper_bound}}]
Consider the class $
\mathcal{L}
=
\left\{
(x,y)\mapsto l(f_j(x),y),
\;
j=1,\ldots,M
\right\}
$. By the standard symmetrization inequality,
\[
\mathbb{E}_{\pi^{\otimes n}}
\Bigl[
R\bigl(\hat f_{\mathrm{ERM}}^{l,\cF}\bigr)\Bigr]
-\min_{1 \leq j \leq M}R(f_j)
\leq
4\,\mathfrak{R}(\mathcal{L}),
\]
where the 
Rademacher complexity
$\mathfrak{R}(\mathcal{L})=\frac1n\,
\mathbb{E}_{\pi^{\otimes n}}
\mathbb{E}_{\varepsilon}
\Bigl[
\sup_{1 \leq j \leq M}
\sum_{i=1}^n \varepsilon_i\,l(f_j(x_i),y_i)
\Bigr]$ and $\varepsilon_i$ are i.i.d. Rademacher random variables.

Since the loss is bounded on an interval of length $h$,
we may subtract the midpoint of this interval from every loss value. This evidently does not affect the expected Rademacher complexity. Consequently, each vector in ${\cal L}$ has entries bounded by $\frac{h}{2}$ and therefore, the Euclidean norm at most $\frac{h}{2}\sqrt{n}$.
Applying Massart's lemma yields
\[
\mathfrak{R}(\mathcal{L})
\le
\frac{h}{2}\sqrt{\frac{2\ln(M)}{n}}=h\sqrt{\frac{\ln(M)}{2n}} .
\]
and, therefore,
\[
\mathbb{E}_{\pi^{\otimes n}}
\Bigl[
R\bigl(\hat f_{\mathrm{ERM}}^{l,\cF}\bigr)\Bigr]
-\min_{1 \leq j \leq M}R(f_j)
\leq
4h\sqrt{\frac{\ln(M)}{2n}}
\]
\end{proof}

\begin{proof}[\textbf{Proof of Theorem~\ref{theorem:erm-lower-bound}}]

    The proof extends the argument of \textcite{lecue2007suboptimalitypenalizedempiricalrisk} to the multiclass setting. For simplicity, assume that $K$ is even. Since $|\mathcal{X}| \geq 2^M$, we may choose a set of $2^M$ distinct points $A=\{x_1,\dots, x_{2^M}\} \subseteq \mathcal{X}$. 
    Fix any bijection between $A$ and the hypercube
    $\Sigma=\{-1,1\}^M$. We may therefore identify $A$ with  $\Sigma$  and consider each point as $x=(x_1,\dots,x_M) \in \Sigma$. Define the marginal distribution  $\pi_X$ of $X$ by taking its coordinates $(x_1,\dots,x_M)$ to be independent Rademacher random variables, and define the conditional (multinomial) distribution of $Y|X=x$ by
   \begin{equation} \label{equation:condprob}
    P(Y|X=x) = \left.
      \begin{cases}
        (\frac{1}{K},...,\frac{1}{K}) & ,\; x_1=-1 \\
         (\frac{2e^\beta}{K(e^\beta+e^{-\beta})}, \frac{2e^{-\beta}}{K(e^\beta+e^{-\beta})}, \frac{2e^\beta}{K(e^\beta+e^{-\beta})}, \frac{2e^{-\beta}}{K(e^\beta+e^{-\beta})},...) & ,\; x_1=1
      \end{cases}
      \right. \;\;\in \mathbb{R}^K
    \end{equation}
    for some $0 < \beta < 1$ to be determined later, 
    and a set of $M$ score functions:
    \[
    f_j(x) = \left.
      \begin{cases}
        (0,...,0) &, \;x_j=-1 \\
         (1,-1,1,-1,...) &, \;x_j=1
      \end{cases}
      \right. \;\in \mathbb{R}^K
    \]
        
   Consider a shifted loss $\hat{l}(f(x),y) = l(f(x),y) - \ln(K)$. This shift obviously does not affect either the excess risk or the empirical risk minimizer. Note that
   \begin{equation} \label{equation:hatl}
   \hat{l}(f_j(x),y)=\mathbb{I}\{x_j=1\}\Big(\ln(\frac{e+e^{-1}}{2})-(-1)^{y+1}\Big)
   \end{equation}
  and, therefore, 
   \begin{gather*}
       \hat{l}(f_1(x),y)-\hat{l}(f_j(x),y) 
       = 
       \mathbb{I}\{x_1 \neq x_j\}x_j\Big((-1)^{y+1} - \ln (\frac{e+e^{-1}}{2}) \Big)
    \end{gather*}
  Consequently, by \eqref{equation:condprob} a straightforward calculus implies
        \begin{gather*}
        R(f_j)-R(f_1)=\frac{\tanh(\beta)}{4},\;\;\;j>1, 
    \end{gather*} 
    where the expectations of $\hat{l}(f_j(x),y)$ are taken with respect to the corresponding joint distribution $\pi(x,y)=\pi_X(x)P(y|X=x)$. Consequently, since $\beta > 0$ we get $R(f_1)=\min_{1 \leq j \leq M}R(f_j)$. Therefore,
    \begin{gather}\label{equation:excess-risk-lower-bound}
        \mathbb{E}_{\pi^{\otimes n}}[R(\hat{f}_{ERM}^{l,\cF})]-\min_{1 \leq j \leq M}R(f_j) 
        = 
        \frac{\tanh(\beta)}{4} P_{\pi^{\otimes n}}(\hat{f}_{ERM}^{\hat{l},\cF} \neq f_1) 
        \geq 
        \frac{\beta}{8} P_{\pi^{\otimes n}}(\hat{f}_{ERM}^{\hat{l},\cF} \neq f_1),
    \end{gather}   
    where the last inequality holds since $\tanh(\beta) > \frac{\beta}{2},\;0 \leq \beta \leq 1$. 
    
    Define the empirical risks $\hat{L}^n(f_j)=\frac{1}{n}\sum_{i=1}^n \hat{l}(f_j(x_i),y_i),\;j=1,\ldots,M$, and 
    let $S=\{n\hat{L}^n(f_1):\{(x_i,y_i)\}_{i=1}^n \in (\mathcal{X} \times \mathcal{Y})^n\}$ be the set of all possible values of $n \hat{L}^n(f_1)$. Evidently, $|S|  <\infty$.
    Conditionally on $Y^n$, $\hat{L}^n (f_j)$ depends only on the coordinate vector $(x_{1j},\ldots,x_{nj}) \in \mathbb{R}^n$. Since these vectors are conditionally independent, $\hat{L}^n(f_1),\ldots,\hat{L}^n(f_M)$ are {\em conditionally} independent given $Y^n$. Moreover, by construction, $\hat{L}^n(f_2),\ldots,\hat{L}^n(f_M)$ are {\em conditionally} i.i.d. Therefore,
    \[
    \begin{split}    
            P(\hat{f}_{ERM}^{\hat{l},\cF} = f_1|Y^n)
        & \leq
        P(n\hat{L}^n(f_1) \leq n\hat{L}^n(f_j) \mid Y^n, \;\forall j) \\
        & =
        \sum_{s \in S}P(n\hat{L}^n(f_1)=s \mid Y^n)P(n\hat{L}^n(f_2) \geq s \mid Y^n)^{M-1}
    \end{split}
    \]
    
    For any $\min(S) \leq s_0 \leq \max(S)$ 
    \[
    \begin{split}   
        \sum_{s \in S} & P(n\hat{L}^n(f_1)=s \mid Y^n)P(n\hat{L}^n(f_2) \geq s \mid Y^n)^{M-1} 
        \\
        &=\sum_{s \in S, s \leq s_0} P(n\hat{L}^n(f_1)=s \mid Y^n)P(n\hat{L}^n(f_2) \geq s \mid Y^n)^{M-1} \\       & + 
        \sum_{s \in S, s > s_0}P(n\hat{L}^n(f_1)=s \mid Y^n)P(n\hat{L}^n(f_2) \geq s \mid Y^n)^{M-1} 
       \\
        & \leq P(n\hat{L}^n(f_1)\leq s_0 \mid Y^n) 
        + 
        P(n\hat{L}^n(f_2) \geq s_0 \mid Y^n)^{M-1}
    \end{split}
    \]
 
    Evidently,
    \begin{gather}\label{equation:probability-lower-bound}
        P(\hat{f}_{ERM}^{\hat{l},\cF} \neq f_1|Y^n) 
        \geq
        1 
        -
        P(n\hat{L}^n(f_1)\leq s_0 \mid Y^n) 
        - 
        P(n\hat{L}^n(f_2) \geq s_0 \mid Y^n)^{M-1}
    \end{gather}
    
    Define $w=(-1)^{y+1}$. Applying Bayes theorem, \eqref{equation:hatl} we have
    \begin{gather*}
        \hat{l}(f_1(X),Y)|(Y=y)  
        =
        \begin{cases}
            \ln(\frac{1+e^{-2w}}{2}), & \text{w.p.}\; \frac{1}{1+\frac{1+e^{-2\beta w}}{2}} 
            \\
            0, & \text{w.p.}\; \frac{1}{1+\frac{2}{1+e^{-2\beta w}}}
        \end{cases}
        ,\;\;\;
        \hat{l}(f_2(X),Y)|(Y=y) =
        \begin{cases}
            \ln(\frac{1+e^{-2w}}{2}), & \text{w.p.}\;\frac{1}{2} \\
            0, & \text{w.p.}\;\frac{1}{2}
        \end{cases}
     \end{gather*} 
yielding    
    \begin{gather*}
        \mathbb{E}_{\mathcal{X}}[\hat{l}(f_1(X),Y)|Y=y] 
        = 
        \ln(\frac{1+e^{-2w}}{2})\frac{1}{1+\frac{1+e^{-2\beta w}}{2}}
        \\
        \mathbb{E}_{\mathcal{X}}[\hat{l}(f_2(X),Y)|Y=y]
        =
        \frac{1}{2}\ln(\frac{1+e^{-2w}}{2}) 
    \end{gather*}

    Take $s_0=n\mathbb{E}_{X^n}[\hat{L}^n(f_1)|Y^n]-\alpha$, where $\alpha$ will be determined later.
    By Hoeffding's inequality,
    \begin{gather}\label{equation:alpha-dependent-upper-bound-f_1-empirical-loss}
        P(n\hat{L}^n(f_1)\leq s_0 \mid Y^n) 
        \leq
        \exp{\left(-{\frac{2\alpha ^ 2}{n\ln^2(\frac{1+e^2}{2})}}\right)}
    \end{gather}
    
    In addition, define
   $$
         \beta_i = \mathbb{E}_{\mathcal{X}}[\hat{l}(f_2(x),Y) \mid Y=y_i] 
         - 
         \mathbb{E}_{\mathcal{X}}[\hat{l}(f_1(X),Y) \mid Y=y_i] = \ln \Big(\frac{2}{1+e^{-2w_i}}\Big) \Big(\frac{1}{1+\frac{1+e^{-2\beta w_i}}{2}}-\frac{1}{2}\Big) 
        \leq 
        \beta
    $$
    for any $0 < \beta < 1$ and
    $$         
    z_i = \mathbb{E}_{\mathcal{X}}[\hat{l}(f_2(X),Y)|Y=y_i]-\hat{l}(f_2(x_i),y_i)
    $$ 
    Note that
    $$
    \mathbb{E}_{\mathcal{X}}[z_i^2 \mid Y=y_i]=\frac{1}{4}\ln^2\Big(\frac{2}{1+e^{-2w_i}}\Big) ,  \;\;\;
        \mathbb{E}_{\mathcal{X}}[|z_i|^3 \mid Y=y_i] = \frac{1}{8}\Big|\ln^3\Big(\frac{2}{1+e^{-2w_i}}\Big)\Big|
    $$
    
    We first lower bound the tail probability of the empirical loss of $f_2$.

    $$
        P\Big(n\hat{L}^n(f_2) \geq s_0 \mid Y^n\Big)
        =
        P\Big(\sum_{i=1}^{n}z_i \leq \alpha + \sum_{i=1}^{n}\beta_i \mid Y^n\Big) 
        \leq 
        P\Big(\frac{\sum_{i=1}^{n}z_i}{\sqrt{\sum_{i=1}^{n}E[z_i^2|y_i]}} \leq \frac{\alpha + \beta n}{c_1\sqrt{n}} \mid Y^n\Big)
    $$
    for some $c_1 > 0$. Plugging in $\alpha = \frac{c_1\sqrt{n\ln(M)}}{2}, \beta=\frac{c_1}{2}\sqrt{\frac{\ln(M)}{n}}$  we have:
    \[   
        P\left(\frac{\sum_{i=1}^{n}z_i}{\sqrt{\sum_{i=1}^{n}E[z_i^2|y_i]}} 
    \leq 
    \frac{\alpha + \beta n}{c_1\sqrt{n}} \mid Y^n\right) 
    = 
    P\Big(\frac{\sum_{i=1}^{n}z_i}{\sqrt{\sum_{i=1}^{n}E[z_i^2|y_i]}} \leq \sqrt{\ln(M)} \mid Y^n\Big)
    \]
    Since $\frac{\mathbb{E}_{\mathcal{X}}[|z_i|^3 \mid Y=y_i]}{\mathbb{E}_{\mathcal{X}}[z_i^2 \mid Y=y_i]} =\frac{1}{2}|\ln \Big(\frac{1 + e^{-2w_i}}{2} \Big)| \leq 1$, using Berry-Esseen theorem and Mills ratio implies:
    \begin{gather*}
        P\Big(\frac{\sum_{i=1}^{n}z_i}{\sqrt{\sum_{i=1}^{n}E[z_i^2|y_i]}} \leq \sqrt{\ln(M)} \mid Y^n\Big) 
        \leq 
        \Phi(\sqrt{\ln(M)})+\frac{1}{\sqrt{n}} 
        \leq 
        1-\frac{1}{\sqrt{32 \pi \ln(M)M}}
    \end{gather*}
    Hence, using the inequality $\ln(1-x)\leq -x$,
    \begin{gather}\label{equation:upper-bound-tail-probability-f_2_empirical-loss}
         P(n \hat{L}^n(f_2) \geq s_0 \mid Y^n)^{M-1} 
         \leq 
         \exp\Big(-{\frac{M-1}{\sqrt{32 \pi M\ln(M)}}}\Big)
    \end{gather}
    
    In addition, plugging the chosen $\alpha=\frac{c_1\sqrt{n\ln(M)}}{2}$ into \eqref{equation:alpha-dependent-upper-bound-f_1-empirical-loss} we have
    \begin{gather}\label{equation:explicit-upper-bound-f_1-empirical-loss}
        P(n\hat{L}^n(f_1)\leq s_0 \mid Y^n) 
        \leq 
        M^{-c_2},
    \end{gather}
    where $c_2=\frac{c_1^2}{2\ln^2 \Big(\frac{1+e^2}{2} \Big)}$. Since the bounds in \eqref{equation:upper-bound-tail-probability-f_2_empirical-loss} and \eqref{equation:explicit-upper-bound-f_1-empirical-loss} hold for every $Y^n$, they hold also unconditionally. Combining \eqref{equation:excess-risk-lower-bound}, \eqref{equation:probability-lower-bound}, \eqref{equation:upper-bound-tail-probability-f_2_empirical-loss}, \eqref{equation:explicit-upper-bound-f_1-empirical-loss} with $\beta=\frac{c_1}{2}\sqrt{\frac{\ln(M)}{n}}$ yields
    \begin{equation*}
    \begin{split}
         \mathbb{E}_{\pi^{\otimes n}} & [R(\hat{f}_{ERM}^{l,\cF})]-\min_{1 \leq j \leq M}R(f_j)
         \geq 
         \frac{c_1}{16}\sqrt{\frac{\ln(M)}{n}}\Big(1- M^{-c_2} - \exp(-{\frac{M-1}{\sqrt{32 \pi M\ln(M)}}})\Big)
    \end{split}
    \end{equation*}
    Take $M_0$ to be the smallest integer such that $M^{-c_2} + \exp(-{\frac{M-1}{\sqrt{32 \pi M\ln(M)}}}) < 1$. Then,  for $M \geq M_0$ we have 
    \[
        \mathbb{E}_{\pi^{\otimes n}} [R(\hat{f}_{ERM}^{l,\cF})]-\min_{1 \leq j \leq M}R(f_j)
        \geq
        c\sqrt{\frac{\ln(M)}{n}}
    \]
    for some $c>0$.
\end{proof}

\begin{proof}[\textbf{Proof of Theorem~\ref{theorem:caew-upper-bound}}]
    Recall that the multinomial logistic loss is $2e^{2B}K$-exp concave (Lemma \ref{lemma:logistic-loss-geometry}). Thus, applying Theorem 4.1 in \textcite{Juditsky_2008} with temperature $\beta=2e^{2B}K$ yields
    \begin{gather}\label{equation:caew-upper-bound-small-K}
        \mathbb{E}_{\pi^{\otimes n}}[R(f_{\hat{\theta}_{CAEW}^{l,\cF}})] - \min_{1 \leq j \leq M}R(f_j) \leq 2e^{2B}~\frac{K\ln(M)}{n}
    \end{gather}

    On the other hand, by convexity of the loss, applying Theorem 2.2 of \textcite{Cesa-Bianchi_Lugosi_2006} to the normalized loss $\tilde{l}(f(x),y) \coloneqq \frac{l(f(x),y) - \ln K + 2B}{4B}$ implies that for $\beta=B\sqrt{\frac{2n}{\ln(M)}}$ 
    \begin{gather*}
        \frac{1}{n}\sum_{i=1}^nl(f_{\hat{\theta}_{AEW,i-1}^{l,\cF}}(x_i),y_i)-\min_{1 \leq j \leq M}\frac{1}{n}\sum_{i=1}^n l(f_j(x_i),y_i) \leq 4B\sqrt{\frac{\ln(M)}{2n}},
    \end{gather*}   
    where $\hat{\theta}_{AEW,i-1}^{l,\cF}$ is the result of applying the AEW procedure defined in \eqref{equation:aew} to the first $i-1$ observations. Note that by definition $\frac{1}{n}\sum_{i=1}^n\hat{\theta}_{AEW,i-1}^{l,\cF}=\hat{\theta}_{CAEW}^{l,\cF}$. Repeating the arguments of the Remark in \textcite[2192]{Juditsky_2008} we have:

    \begin{gather}\label{equation:caew-upper-bound-large-K}
        \mathbb{E}_{\pi^{\otimes n}}[R(f_{\hat{\theta}_{CAEW}^{l,\cF}})] - \min_{1 \leq j \leq M}R(f_j) 
        \leq 
        4B\sqrt{\frac{\ln(M)}{2n}}
    \end{gather} 

   Taking $c(B) = \max \{2e^{2B}, 4B \}$, \eqref{equation:caew-upper-bound-small-K} and \eqref{equation:caew-upper-bound-large-K} imply
    \begin{gather*}
        \mathbb{E}_{\pi^{\otimes n}}[R(f_{\hat{\theta}_{CAEW}^{l,\cF}})] - \min_{1 \leq j \leq M}R(f_j) 
        \leq 
        c(B)\min \Big\{\frac{K\ln(M)}{n}, \sqrt{\frac{\ln(M)}{n}}\Big\}
    \end{gather*}

\end{proof}

\begin{proof}[\textbf{Proof of Theorem~\ref{theorem:proper-estimators-lower-bound}}]
    Let $N = \lfloor \log_2(M) \rfloor$. Consider a set of $N+1$ distinct points $\mathcal{X}=\{x_1, ... x_{N+1}\}$ and 
    define the marginal distribution on $\mathcal X$ as follows:
    \[
         P(x_i) = 
            \begin{cases}
                w, & 1 \leq i \leq N \\
                1-wN, & i=N+1
            \end{cases}
    \]
    Consider the vectors $u=(\frac{1}{2},\frac{1}{2},0,\ldots,0)^\top \in \mathbb{R}^K,\;v=(\epsilon,-\epsilon,0,\ldots,0)^\top \in \mathbb{R}^K$ and $z=(a,-a,0,\ldots,0)^\top \in \mathbb{R}^K$, 
    where $w>0,\;a > 0$ and $0 < \epsilon < \frac{1}{2}$ will be defined later.
    
    Define two distributions over $\mathcal{Y}=\{1,..,K\}$ with probabilities
    \[
         Q_1 = u+v, \;Q_{-1} = u-v
    \]
    and the set of distributions $\Pi = \{\pi_{\sigma}:~\sigma \in \Sigma \}$
    over $\mathcal X \times \mathcal Y$ indexed by the hypercube $\Sigma = \{-1,1\}^{N}$ as follows
    
    \[
        \pi_\sigma(x_i, y) =
            \begin{cases}
                w Q_{\sigma_i}(y), &  1 \leq i \leq N \\
                (1-wN)  Q_1(y), &  i=N+1
            \end{cases}
    \]
    Denote by $\pi_{\sigma,x_i}$ the conditional distribution of $y|(X=x_i),\;i=1,\ldots,N+1$ induced by $\pi_{\sigma}$, i.e.,
    \[
        \pi_{\sigma,x_i}(y) =
            \begin{cases}
                Q_{\sigma_i}(y), &  1 \leq i \leq N \\
                Q_1(y), &  i=N+1
            \end{cases}
    \]

    For each $\sigma \in \Sigma$ define a score function
    
    \[
        f_{\sigma}(x_i) = 
            \begin{cases}
                z, &  1\leq i \leq N \text{ and } \sigma_i = 1 \\
                -z, & 1\leq i \leq N \text{ and } \sigma_i = -1 \\
                z, &  i=N+1 \\
            \end{cases}
    \]
    and consider a dictionary of size $M$: $\cF=\{f_{\sigma}:  \sigma \in \Sigma \} \cup \{f_{2^{N}+1}, \dots, f_{M}\}$, where the padding score functions $f_j, \;j=2^{N}+1,\ldots,M$ are any distinct convex combinations of $f_\sigma$ such that $\lVert f_j(x) \rVert_\infty \leq a$. Evidently, $\operatorname{conv}(\cF)=conv\{f_\sigma\}$.

    By the design of the dictionary $\cF$, for any $\hat{\theta} \in \Delta_M$ and the corresponding aggregate $\hat{f}(x)=\sum_{j=1}^{M} \hat{\theta}_j f_j(x) \in \operatorname{conv}(\cF)$ we have
    \[
    \hat{f}(x_i) = (\lambda_{\hat{\theta} i},-\lambda_{\hat{\theta} i},0,\ldots,0)^T,
    \]
    where $|\lambda_{\hat{\theta} i}| \leq a$.

    By Lemma \ref{lemma:loss-strongly-convex-and-smooth-wrt-function}, for every $\sigma \in \Sigma$, every $\hat{f} \in \operatorname{conv}(\cF)$  and any $x,y$:
    \begin{gather}\label{equation:taylor-lower-bound-general-case}
        l(\hat{f}(x),y)-l(f_{\sigma}(x),y) 
        \geq 
        \nabla l(f_{\sigma}(x),y)^T(\hat{f}(x)-f_{\sigma}(x)) + \frac{e^{-2a}}{2K} \lVert \hat{f}(x)-f_{\sigma}(x) \rVert_2^2,
    \end{gather}
    where $\nabla l(f_{\sigma}(x),y) =-e_{y}+\operatorname{softmax}(f_{\sigma}(x))$.
    
    Hence, 
  \begin{equation} \label{equation:conditional_risk}
        \mathbb{E}_{\pi_{\sigma,x}}[l(\hat{f}(x)) 
        - l(f_{\sigma}(x))] 
        \geq 
        \left(-\pi_{\sigma, x} + \operatorname{softmax}(f_{\sigma}(x))\right)^T(\hat{f}(x)-f_{\sigma}(x)) 
        + 
        \frac{e^{-2a}}{2K}\lVert \hat{f}(x)-f_{\sigma}(x) \rVert_2^2
    \end{equation}
    
    Plugging the values a straightforward calculus implies
    
    \begin{equation} \label{equation:gradient-general-form}
         \Big(-\pi_{\sigma,x_i} + \operatorname{softmax}(f_{\sigma}(x_i))\Big)^T(\hat{f}(x_i)-f_{\sigma}(x_i)) 
         =
         \Big(2\epsilon - \frac{e^a-e^{-a}}{e^{a} + e^{-a}+K-2}\Big)\Big(a -  \sigma_i \lambda_{{\hat \theta} i}\Big) 
    \end{equation}
       
    For $1 \leq i \leq N$. Notice that $a  \geq \sigma_i \lambda_{{\hat \theta} i}$ since $|\lambda_{{\hat \theta} i}| \leq a$.
    
    Consider two cases:
    \newline \textbf{Case 1}. $K^2 \leq \frac{32n}{\log_2(M)}$.
    
    Taking $\epsilon = \frac{1}{K}$ and $a=\ln(2)$ in \eqref{equation:gradient-general-form} yields
    \begin{equation}\label{equation:small-K-linear-term-non-negative}
        \left(-\pi_{\sigma,x_i} + \operatorname{softmax}(f_\sigma(x_i))\right)^T(\hat{f}(x_i)-f_\sigma(x_i)) 
        \geq
        0
    \end{equation}
        for any $1 \leq i \leq N$. Taking the expectation with respect to $x$,  \eqref{equation:taylor-lower-bound-general-case} and \eqref{equation:small-K-linear-term-non-negative} imply that for any $\hat{f} \in \operatorname{conv}(\cF)$:
        \[
    R(\hat{f}) - R(f_{\sigma}) \geq \frac{1}{8K}\mathbb{E}_{\pi_{\sigma}}\lVert \hat{f}(x)-f_{\sigma}(x) \rVert_2^2
    \]
    Define the metric $d(f,g) =\sqrt{ \sum_{i=1}^{N+1} \sum_{k=1}^{K}P_x(x_i) \frac{1}{K} (f_k(x_i) - g_k(x_i))^2}$.  We have
    \begin{gather*}
        \inf_{\hat{f} \in \operatorname{conv}(\cF)}\sup_{\pi} \mathbb{E}_{\pi^{\otimes n}}\Big[ R(\hat{f}) 
        - 
        \min_{1 \leq j \leq M}R(f_j)\Big] 
        \geq 
       \frac{1}{8}\inf_{\hat{f} \in \operatorname{conv}(\cF)}\sup_{\pi_{\sigma} \in \Pi}\mathbb{E}_{\pi_{\sigma}^{\otimes n}}\Big[d^2(f_{\sigma}, \hat{f})\Big]
    \end{gather*}
            
    For any $\sigma, \tilde{\sigma}$:
    \begin{gather}\label{equation:squared-distance-small-K}
        d^2(f_{\sigma},f_{\tilde{\sigma}}) 
        = 
       \frac{1}{K}\mathbb{E}_{\pi_{\sigma}}|\lVert_{\tilde{\sigma}}(x)-f_{\sigma}(x) \rVert_2^2
        =
        \rho(\sigma, \tilde{\sigma}) \frac{w}{K}\lVert 2z \rVert_2^2
        =
        \rho(\sigma, \tilde{\sigma}) \frac{w}{K} 8 \ln^2(2),
    \end{gather}
    where $\rho(\sigma, \tilde{\sigma})$ is the Hamming distance. 
    In addition, for any $\sigma, \tilde{\sigma}$ such that $\rho(\sigma, \tilde{\sigma}) = 1$, the squared Hellinger distance $H^2(\pi_{\sigma}, \pi_{\tilde{\sigma}})$ between $\pi_\sigma$ and $\pi_{\tilde \sigma}$ satisfies 
        \begin{gather*}
        H^2(\pi_{\sigma}, \pi_{\tilde{\sigma}}) = \sum_{i=1}^{N+1} \sum_{k=1}^{K}(\sqrt{\pi_{\sigma}(x_i, k)} - \sqrt{\pi_{\tilde{\sigma}}(x_i, k)})^2 = 2\Big(\sqrt{w(\frac{1}{2} + \frac{1}{K})}  - \sqrt{w(\frac{1}{2} - \frac{1}{K}})\Big)^2 = \\
        2w\Big(1 - 2\sqrt{\frac{1}{4} - \frac{1}{K^2}}\Big) \leq  2w\Big(1-(1- \frac{4}{K^2})\Big) = \frac{8w}{K^2}
    \end{gather*}
    Taking $w = \frac{K^2}{32n}$ yields
    \begin{gather*}
        H^2(\pi_{\sigma}^{\otimes{n}}, \pi_{\tilde{\sigma}}^{\otimes{n}})= 
2\left(
1-
\left(
1-\frac{H^2(\pi_\sigma,\pi_{\tilde \sigma})}{2}
\right)^n
\right)
        \leq
        2\Big(1-(1-\frac{1}{8n})^n\Big) 
        \leq 
        \frac{1}{4}
    \end{gather*}
    Using the known inequality to lower bound the total variation affinity $\lVert P \wedge Q \rVert \geq 1-H(P,Q)$ we get
    \begin{gather}\label{equation:affinity-small-K}
        \min_{
            \sigma_1, \sigma_2:\rho(\sigma_1, \sigma_2)=1
        } 
        \lVert 
            \pi_{\sigma_1}^{\otimes n} \wedge \pi_{\sigma_2}^{\otimes n}
        \rVert
        \geq  
        \frac{1}{2}
    \end{gather}

    Applying Lemma 4 of \textcite{Cai_2010}, a variant of Assouad's lemma, together with \eqref{equation:squared-distance-small-K} and \eqref{equation:affinity-small-K}, yields
    \begin{gather}\label{equation:squared-distance-lower-bound}
        \inf_{\hat{f} \in \operatorname{conv}(\cF)}\sup_{\sigma} \mathbb{E}_{\pi_{\sigma}^{\otimes n}}[d^2(\hat{f} , f_{\sigma})] 
        \geq 
        \frac{1}{4}\frac{K\ln(2)^2}{4n}\frac{\lfloor \log_2(M) \rfloor}{2}\frac{1}{2} \geq \frac{\ln^2(2)K\log_2(M)}{128n}
    \end{gather}

    \textbf{Case 2}. $K^2 \geq \frac{9n}{\log_2(M)}$.         
    For this case we use the linear part of the Taylor expansion in \eqref{equation:conditional_risk} to get the lower bound for the the excess risk. From \eqref{equation:conditional_risk},    
    \begin{gather*}
        \mathbb{E}_{\pi_{\sigma,x}}[l(\hat{f}(x))] - \mathbb{E}_{\pi_{\sigma,x}}[l(f_{\sigma}(x))] \geq
        \Big(-\pi_{\sigma,x} + \operatorname{softmax}(f_{\sigma}(x))\Big)^T(\hat{f}(x)-f_{\sigma}(x))
    \end{gather*}
    
    Taking $\epsilon = \sqrt{\frac{\log_2(M)}{n}}, a=1$ and plugging them into \eqref{equation:gradient-general-form} for $K \geq 3\sqrt{\frac{n}{\log_2(M)}}$ yields
    \[
    \begin{split}
         \Big(& -\pi_{\sigma,x_i} +\operatorname{softmax}(f_{\sigma}(x_i))\Big)^T
         (\hat{f}(x_i) - f_{\sigma}(x_i)) 
         \geq 
         (1-\sigma_i \lambda_{{\hat \theta}i}) \sqrt{\frac{\log_2(M)}{n}} \\
        &=
        |\sigma_i - \lambda_{{\hat \theta}i}|\sqrt{\frac{\log_2(M)}{n}},
    \end{split}
    \]
  where the last equality holds since $\sigma_i=\pm 1$ and $| \lambda_{{\hat \theta}i}| \leq 1$.  

    Therefore, 
    \begin{gather}\label{equation:expected-risk-l1-lower-bound}
       \inf_{\hat{f} \in \operatorname{conv}(\cF)} \sup_{\pi} \mathbb{E}_{\pi^{\otimes n}}\Big[ R(\hat{f}) 
        - 
        \min_{f \in \cF}R(f)\Big] 
        \geq 
        w \sqrt{\frac{\log_2(M)}{n}} \inf_{\hat{\lambda}} \sup_{\pi_{\sigma} \in \Pi}  \mathbb{E}_{\pi_{\sigma}^{\otimes n}} \Big[\sum_{i=1}^{N}|\sigma_i - \lambda_{{\hat \theta}i}| \Big]
    \end{gather}
    In addition, for any $\sigma, \tilde{\sigma}$ such that $\rho(\sigma, \tilde{\sigma}) = 1$, 
    \[
    \begin{split}
        H^2(\pi_{\sigma}, \pi_{\tilde{\sigma}}) 
        & = 
        \sum_{i=1}^{N+1} \sum_{k=1}^{K}(\sqrt{\pi_{\sigma}(x_i, k)} - \sqrt{\pi_{\tilde{\sigma}}(x_i, k)})^2 
        = 
        2\Big(\sqrt{w(\frac{1}{2} + \sqrt{\frac{\log_2(M)}{n}})}  - \sqrt{w(\frac{1}{2} - \sqrt{\frac{\log_2(M)}{n}}})\Big)^2  \\
        & = 
        2w\Big(1 - 2\sqrt{\frac{1}{4} - \frac{\log_2(M)}{n}}\Big) 
        \leq  
        2w\Big(1-(1- \frac{4\log_2(M)}{n})\Big) 
        = 
        \frac{8w\log_2(M)}{n}
    \end{split}
    \]
    Hence, taking $w = \frac{1}{16\log_2(M)}$ we have
    \begin{gather*}
        H^2(\pi_{\sigma}^{\otimes{n}}, \pi_{\tilde{\sigma}}^{\otimes{n}}) 
        \leq 
        2\left(1-(1-\frac{1}{4n})^n\right) 
        \leq 
        \frac{1}{2}
    \end{gather*}
    
    Applying Lemma 1 of \textcite{lecue2007suboptimalitypenalizedempiricalrisk}, another version of Assouad's lemma, to $\omega_i=\frac{\sigma_i + 1}{2}, \tilde{\omega}_i=\frac{\lambda_{{\hat \theta}i}+1}{2}$
    \begin{gather}\label{equation:l1-lower-bound}
        \inf_{\lambda_{\hat \theta}}\sup_{\sigma \in \Sigma} \mathbb{E}_{\pi_{\sigma}^{\otimes n}}\sum_{i=1}^{N}|\sigma_i - \lambda_{{\hat \theta}i}| 
        \geq 
        \frac{9N}{32} 
        \geq
        \frac{9\log_2(M)}{64}
    \end{gather}
    which together with \eqref{equation:expected-risk-l1-lower-bound} implies
    
    \begin{gather} \label{equation:lower_bound_case2}
        \inf_{\hat{f} \in \operatorname{conv}(\cF)} \sup_{\sigma \in \Sigma}\mathbb{E}_{\pi_{\sigma}^{\otimes{n}}}\Big[R(\hat{f}) - R(f_{\sigma})\Big] 
        \geq  
        \frac{1}{120}\sqrt{\frac{\log_2(M)}{n}}
    \end{gather}
    
    Combining \eqref{equation:squared-distance-lower-bound} with \eqref{equation:lower_bound_case2} completes the proof.

\end{proof}

\begin{proof}[\textbf{Proof of Theorem~\ref{theorem:optimal-ms-estimator}}] 
    Consider the loss $l(f_\theta^P(x),y)$. For any $\tilde{\theta} \in \Delta_M$, define
    $$
    Q_1(\theta,\tilde{\theta},x,y)=l(f_{\theta}^P(x),y)-l(f_{\tilde \theta}^P(x),y)=\ln\left( \frac{\sum_{j=1}^M\tilde{\theta}_jp_{jy}(x)}{\sum_{j=1}^M\theta_jp_{jy}(x)}\right)
    $$
    Evidently, as a function of $\theta$,  $Q_1(\theta,\tilde{\theta},x,y)$ is 1-exp-concave on $\Delta_M$ for all $x, y$ and $\tilde{\theta}$. Hence, applying Theorem 4.1 in \textcite{Juditsky_2008} for $\hat{f}^P_{CAEW}$ in \eqref{equation:probability_aggregate} with $\beta=1$ yields
    $$
    \mathbb{E}_{\pi^{\otimes n}}[R(\hat{f}^P_{CAEW})] \leq  \min_{1 \leq j \leq M}R(f_j) + \frac{\ln(M)}{n}
    $$

\end{proof}

\begin{proof}[\textbf{Proof of Theorem~\ref{theorem:lower-bound-ms-general-case}}]
The proof follows the lines of the proof of
Theorem~\ref{theorem:proper-estimators-lower-bound}
and is based on an application of Lemma~4 of \textcite{Cai_2010},
a version of Assouad's lemma, although the construction of the
distributions and score functions is somewhat different.
    
    We assume again for simplicity that $K$ is even.
    Define the vector $g=(-1,1,\ldots,-1,1)^\top \in \mathbb{R}^K$ and let  $Q_1 = \operatorname{softmax}(g), \;
        Q_{-1} = \operatorname{softmax}(-g)$.
    
    Consider distinct points $x_1, ..., x_{N+1}$ where $N = \lfloor\log_2(M)\rfloor$, and define the marginal distribution:
    
    \begin{gather*}
        P_X(x_i) =
            \begin{cases}
                w, & 1 \leq i \leq N \\
                1-wN, & i = N + 1
            \end{cases}
    \end{gather*}
    where $w > 0$ will be determined later. For every $\sigma \in \Sigma=\{-1, 1\}^{N}$  define a score function
    \begin{gather*}
        f_{\sigma}(x_i) =
            \begin{cases}
                \sigma_{i}g, & 1 \leq i \leq N \\
                g & i = N + 1
            \end{cases}
    \end{gather*}
    Hence, $Q_{\sigma_ i}=\operatorname{softmax}(f_\sigma(x_i))$. 
The corresponding joint distribution is
\[
\pi_\sigma(y,x_i)
=
P_X(x_i)\,
\operatorname{softmax}(f_\sigma(x_i))_y,
\]
while the conditional distribution is
\begin{equation}
\label{equation:conditional_distribution}
\pi_{\sigma,x_i}
=
P_{\pi_\sigma}\!\left(\,\cdot\,\middle|X=x_i\right)
=
\operatorname{softmax}(f_\sigma(x_i)).
\end{equation}
 
Consider a dictionary of size $M$: $\cF=\{f_{\sigma}:  \sigma \in \Sigma \} \cup \{f_{2^{N}+1}, \dots, f_{M}\}$, where the distinct padding score functions $f_j, \;j=2^{N}+1,\ldots,M$ satisfy $\lVert f_j(x) \rVert_\infty \leq 1$. 
  
   For a given $\pi_\sigma$, for every score function $f$, \eqref{equation:conditional_distribution} implies
\[
R(f)-R(f_\sigma)
=
\mathbb{E}_{X}\!\left[
\operatorname{KL}\!\left(
\pi_{\sigma,x}\,\parallel \,\operatorname{softmax}(f(x))
\right)
\right]
\geq 0,
\]
and therefore
\[
\min_{f\in \cF}R(f)=R(f_\sigma).
\]

    Define a set of distributions $\Pi = \{\pi_{\sigma}: \sigma \in \Sigma\}$. We have
        \begin{gather*}
        \inf_{\hat{f}}\sup_{\pi} \Big[\mathbb{E}_{\pi^{\otimes n}}[ R(\hat{f})] 
        - 
        \min_{f \in \cF}R(f)\Big] 
        \geq 
        \inf_{\hat{f}}\sup_{\pi_{\sigma} \in \Pi} \Big[\mathbb{E}_{\pi_{\sigma}^{\otimes n}}[R(\hat{f})] - \min_{f \in \cF}R(f)\Big]
        = \\
       \inf_{\hat{f}}\sup_{\pi_{\sigma} \in\Pi}\mathbb{E}_{\pi_{\sigma}^{\otimes n}}\Big[\mathbb{E}_X[KL(\pi_{\sigma,x} \parallel \pi_{\hat{f},x})]\Big]
       \geq 
       \inf_{\hat{f}}\sup_{\pi_{\sigma} \in \Pi}\mathbb{E}_{\pi_{\sigma}^{\otimes n}}\Big[H^2(\pi_{\sigma}, \pi_{\hat{f}})\Big]
    \end{gather*}
        By construction, a straightforward calculus implies that for every $\sigma_1, \sigma_2 \in \Sigma$
    \begin{equation} \label{equation_hellinger}
    \begin{split}
        H^2(\pi_{\sigma_1}, \pi_{\sigma_2}) 
        & = \sum_{i=1}^{N+1} \sum_{k=1}^{K}(\sqrt{\pi_{\sigma_1}(x_i, k)} - \sqrt{\pi_{\sigma_2}(x_i, k)})^2 \\
        & = 
        w\rho(\sigma_1, \sigma_2)\left(2-\frac{4e}{e^2+1}\right) 
        \leq 
         w\rho(\sigma_1, \sigma_2)
    \end{split}
    \end{equation}
    Take $w=\frac{1}{4n}$. Then for
    $\rho(\sigma_1, \sigma_2)=  1$ we have
    \begin{gather*}
         H^2(\pi_{\sigma_1}^{\otimes n}, \pi_{\sigma_2}^{\otimes n}) 
         = 
         2\left(1-(1-\frac{H^2(\pi_{\sigma_1}, \pi_{\sigma_2})}{2})^n\right)
         \leq 2\left(1-(1-\frac{1}{8n})^n\right) 
         \leq 
         \frac{1}{4}
    \end{gather*}
    
    Hence, the total variation affinity
    \begin{gather}\label{equation:total-variation-lpwer-bound}
        \min_{\sigma_1, \sigma_2:\rho(\sigma_1,\sigma_2)=1}\lVert \pi_{\sigma_1}^{\otimes n} \wedge \pi_{\sigma_2}^{\otimes n} \rVert
        \geq
        1 - \max_{\sigma_1, \sigma_2:\rho(\sigma_1,\sigma_2)=1} H(\pi_{\sigma_1}^{\otimes n}, \pi_{\sigma_2}^{\otimes n})
        \geq  
        \frac{1}{2}
    \end{gather}

    On the other hand, from \eqref{equation_hellinger}, 
    \begin{gather}\label{equation:hellinger-distance-lower-bound}
        H^2(\pi_{\sigma_1}, \pi_{\sigma_2}) 
        = 
        \frac{1}{4n}\rho(\sigma_1, \sigma_2)\left(2-\frac{4e}{e^2+1}\right) 
        \geq 
        \frac{\rho(\sigma_1, \sigma_2)}{8n}
    \end{gather}
    Applying Lemma 4 of \textcite{Cai_2010}, 
    \eqref{equation:total-variation-lpwer-bound} and \eqref{equation:hellinger-distance-lower-bound} yield
        \begin{gather*}
        \inf_{\hat{f}}\sup_{\pi} \mathbb{E}_{\pi^{\otimes n}}\Big[ R(\hat{f}) 
        - 
        \min_{f \in \cF}R(f)\Big] 
        \geq 
        \inf_{\hat{f}}\sup_{\pi_{\sigma} \in\Pi}\mathbb{E}_{\pi_{\sigma}^{\otimes n}}\Big[H^2(\pi_{\sigma}, \pi_{\hat{f}})\Big]
        \geq 
        \frac{\log_2(M)}{256n}
    \end{gather*}
\end{proof}

\begin{proof}[\textbf{Proof of Theorem~\ref{theorem:convex-erm-0/1}}]
Consider the class
\[
\mathcal H=
\left\{
x\mapsto
\operatorname*{arg\,max}_{1\le k\le K}f_{\theta,k}(x):
\theta\in\Delta_M
\right\}  \subset \left\{
x\mapsto
\operatorname*{arg\,max}_{1\le k\le K}f_{\theta, k}(x):
\theta \in \mathbb{R}^M\right\} =:\tilde{\mathcal{H}}
\]
Applying Theorem~29.7 of \textcite{Shalev-Shwartz_Ben-David_2014} for $\Psi(x,k)=\left(f_{1k}(x),\ldots,f_{Mk}(x)\right) \in \mathbb{R}^M$, implies that the Natarajan dimension
\[
\operatorname{Ndim}(\mathcal H)\le \operatorname{Ndim}(\tilde{\mathcal H}) \leq M
\]
From Theorem~6 of \textcite{JMLR:v16:daniely15a} it then follows that, with probability at least \(1-\delta\),
\[
R_{0/1}\!\left(\hat f^{\,l_{0/1},\operatorname{conv}(\cF)}_{\mathrm{ERM}}\right)
-
\inf_{f\in\operatorname{conv}(\cF)}R_{0/1}(f)
\le
c\sqrt{\frac{M\ln K+\ln(1/\delta)}{n}}
\]
for some universal constant \(c>0\). The bound in expectation follows by integrating the above tail probability. 

\end{proof}

\begin{proof}[\textbf{Proof of Theorem~\ref{theorem:erm-convex-hull-upper-bound}}]
    By Lemma \ref{lemma:logistic-loss-geometry}, for all $\theta, \tilde{\theta} \in \Delta_M$ we have
    $$
    \frac{1}{n}\sum_{i=1}^{n}\Big(l(f_{\theta}(x_i),y_i) - l(f_{\tilde{\theta}}(x_i),y_i) \Big)^2 \leq \frac{2}{n}\sum_{i=1}^n (\theta-\tilde{\theta})^TF(x_i)^TF(x_i)(\theta-\tilde{\theta})
    $$
    and
   $$ 
        L^n(f_\theta) - L^n(f_{\tilde \theta}) \geq 
        \nabla L^n(\tilde{\theta})^T(\theta-\tilde{\theta})+\frac{1}{2e^{2B}K}(\theta-\tilde{\theta})^T \frac{1}{n} \sum_{i=1}^nF(x_i)^TF(x_i) (\theta-\tilde{\theta})
    $$ 
    Hence, the empirical multinomial logistic risk $L^n$ satisfies Assumption~3.1 in \textcite{puchkin2023exploringlocalnormsexpconcave} with respect to the semi-norm induced by the matrix $\frac{1}{n}\sum_{i=1}^{n}F(x_i)^TF(x_i)$. Theorem 3.4 in the same paper implies then that for every $\delta > 0$, with probability of at least $1-\delta$
    \[
    R(\hat{f}_{ERM}^{l,\operatorname{conv}(\cF)}) - \inf_{f \in \operatorname{conv}(\cF)}R(f) \lesssim  \frac{2e^{2B}K(M + \ln(\frac{1}{\delta}))}{n}
    \]
    Integrating the tail probability yields the bound in expectation:
    \begin{equation} \label{equation:firstterm}
    \mathbb{E}_{\pi^{\otimes n}}[R(\hat{f}_{ERM}^{l,\operatorname{conv}(\cF)})] - \inf_{f \in \operatorname{conv}(\cF)}R(f) \leq  c_1(B) \frac{KM}{n}
    \end{equation}

On the other hand,
let
\[
\mathcal L_{\Delta_M}
=
\Bigl\{(x,y)\mapsto \ell(F(x)\theta,y):\theta\in\Delta_M\Bigr\},
\qquad
\ell(f,y)=-f_y+\ln\!\Bigl(\sum_{k=1}^K e^{f_k}\Bigr).
\]
By the standard symmetrization inequality,
\[
\mathbb{E}_{\pi^{\otimes n}}
\Bigl[
R\bigl(\hat f_{\mathrm{ERM}}^{\,l,\mathrm{conv}(\cF)}\bigr)\Bigr]
-\inf_{f\in \mathrm{conv}(\cF)}R(f)
\le
4\,\mathfrak R(\mathcal L_{\Delta_M}),
\]
where the Rademacher complexity 
\begin{equation*}
\begin{split}
\mathfrak R(\mathcal L_{\Delta_M})
& =
\frac1n\,
\mathbb{E}_{\pi^{\otimes n}}
\mathbb{E}_{\varepsilon}
\Bigl[
\sup_{\theta\in\Delta_M}
\sum_{i=1}^n \varepsilon_i\,\ell(F(x_i)\theta,y_i)
\Bigr] \\
& \leq \frac1n
\mathbb{E}_{\pi^{\otimes n}}\mathbb{E}_\varepsilon
\sup_{\theta\in\Delta_M}
\sum_{i=1}^n \varepsilon_i\bigl(-f_{\theta,y_i}(x_i)\bigr) +
\frac1n
\mathbb{E}_{\pi^{\otimes n}}\mathbb{E}_\varepsilon
\sup_{\theta\in\Delta_M}
\sum_{i=1}^n \varepsilon_i\,\ln\!\Bigl(\sum_{k=1}^K e^{(F(x_i)\theta)_k}\Bigr) \\
&=: T_1 + T_2
\end{split}
\end{equation*}
and $\varepsilon_i$ are i.i.d. Rademacher random variables.

For the first term, since \(\Delta_M\) is the simplex, the supremum over \(\theta\) is attained at a vertex: 
\[
T_1=\frac1n
\mathbb{E}_{\pi^{\otimes n}}\mathbb{E}_\varepsilon
\max_{1\le j\le M}
\sum_{i=1}^n \varepsilon_i\bigl(-f_{j,y_i}(x_i)\bigr).
\]
Because \(|f_{j,k}(x)|\le B\), Massart's lemma implies
\[
T_1 \le B\sqrt{\frac{2\ln M}{n}}.
\]

For the second term, one can easily verify that the log-sum-exp function is
$\frac{e^{2B}}{\sqrt K}$-Lipschitz with respect to the Euclidean norm on $[-B,B]^K$.
Therefore, by vector contraction inequality \parencite[Corollary 4]{maurer2016vectorcontractioninequalityrademachercomplexities} ,
\[
T_2
\le
\frac{\sqrt2\,e^{2B}}{n\sqrt K}\,
\mathbb{E}_{\pi^{\otimes n}}\mathbb{E}_\varepsilon
\sup_{\theta\in\Delta_M}
\sum_{i=1}^n\sum_{k=1}^K
\varepsilon_{ik}\,f_{\theta,k}(x_i),
\]
where \(\varepsilon_{ik}\) are i.i.d. Rademacher variables. Again, since \(\Delta_M\) is a simplex,
\[
\sup_{\theta\in\Delta_M}
\sum_{i=1}^n\sum_{k=1}^K
\varepsilon_{ik}\,f_{\theta,k}(x_i)
=
\max_{1\le j\le M}
\sum_{i=1}^n\sum_{k=1}^K
\varepsilon_{ik}\,f_{jk}(x_i).
\]
and Massart's lemma yields
\[
\mathbb{E}_{\pi^{\otimes n}}\mathbb{E}_\varepsilon
\max_{1\le j\le M}
\sum_{i=1}^n\sum_{k=1}^K
\varepsilon_{ik}\,f_{jk}(x_i)
\le
B\sqrt{2nK\ln M}.
\]

Consequently,
\[
T_2
\le
2Be^{2B}\sqrt{\frac{\ln M}{n}}.
\]
and, therefore,
\[
\mathfrak R(\mathcal L_{\Delta_M}) \leq 
c_2(B)\sqrt{\frac{\ln M}{n}}
\]
for a constant \(c_2(B)>0\) depending only on \(B\). Hence,
\begin{equation}  \label{equation:secondterm}
\mathbb{E}_{\pi^{\otimes n}}
\Bigl[
R\bigl(\hat f_{\mathrm{ERM}}^{\,l,\mathrm{conv}(\cF)}\bigr)\Bigr] -
\inf_{f\in \mathrm{conv}(\cF)}R(f)
\leq 4c_2(B)\sqrt{\frac{\ln M}{n}}.
\end{equation}

Combining \eqref{equation:firstterm} and \eqref{equation:secondterm} completes the proof.

\end{proof}

\begin{proof}[\textbf{Proof of Theorem~\ref{theorem:convex-aggregation-lower-bound}}]
    Assume for simplicity that $K$ is even. 
    Since $|\mathcal{X}| \geq 2^M$, we may choose $2^M$ distinct points $T \coloneqq \{x_1,\dots, x_{2^M}\} \subseteq \mathcal{X}$. Fix any bijection between $T$ and the hypercube
    $\Sigma=\{-1,1\}^M$. We may therefore identify $T$ with  $\Sigma$  and consider each point as $x=(x_1,\dots,x_M) \in \Sigma$. Assume the marginal uniform distribution $\pi_X(x)=\frac{1}{2^M},\; x \in \Sigma$.

    For any function $f: \mathcal{X} \to \mathbb{R}^K$, let $p_f(x) =\operatorname{softmax}(f(x)) \in \Delta_K$. 
    Define the score functions:
    \begin{gather*}
        f_j(x)=\epsilon x_j v
    \end{gather*}
    for some $0 < \epsilon < 1$ to be determined later, where  the vector $v=(1,-1,1,-1, \dots)^T \in \mathbb{R}^K$. 
    
    For any $\theta \in \Delta_M$ consider $f_\theta(x)=\sum_{j=1}\theta_j f_j(x)=\epsilon \theta^\top x v$. To simplify notations we denote $p_{\theta}(x)=p_{f_{\theta}}(x)$.
       
    Thus,           
       \begin{equation}   \label{equation:p_thetak}   
        p_{\theta,k}(x)
=
\begin{cases}
\frac{2}{K}
\frac{e^{\epsilon\theta^\top x}}
{e^{\epsilon\theta^\top x}+e^{-\epsilon\theta^\top x}},
& k\;{\rm is\;odd}, \\
\frac{2}{K}
\frac{e^{-\epsilon\theta^\top x}}
{e^{\epsilon\theta^\top x}+e^{-\epsilon\theta^\top x}},
& k\;{\rm is\; even}.
\end{cases}
\end{equation}

    For every $1 \leq m \leq \frac{M}{2}$, consider a set of $m$-sparse vectors 
    $
    \Theta^M_m := \{\theta \in \{0,\frac{1}{m}\}^M:  \lVert \theta \rVert_0 = m\}
    $
    By Lemma 8.3 of \textcite{rigollet2010exponentialscreeningoptimalrates}, there exists a subset $\tilde{\Theta}_{m}^{M} \subseteq  \Theta^M_m$ such that for any $\theta_1, \theta_2 \in \tilde{\Theta}_{m}^{M}$
    \begin{gather}\label{equation:hamming-distance-lower-bound-theorem-9}
        \rho(\theta_1, \theta_2) 
        \geq 
        \frac{m}{4} 
       \end{gather}
    and
    \begin{gather}\label{equation:set-size-lower-bound-theorem-9}
        \ln|{\tilde{\Theta}_{m}^{M}}| 
        \geq 
        c_1m \cdot \ln(1+M/m)
    \end{gather}
    for some $0<c_1<1$.

    Let $\cF_{\tilde{\Theta}_{m}^{M}}  = \{f_\theta: \theta \in {\tilde{\Theta}_{m}^{M}}  \}$ and $P_{\tilde{\Theta}_{m}^{M}} =\{p_\theta:  \theta \in {\tilde{\Theta}_{m}^{M}}  \}$. 
    Define the distribution $\pi_{\theta}(x,k)=\pi_{X}(x) p_{\theta, k}(x)$ on $\Sigma \times \mathcal{Y}$.
        Note that under $\pi_{\theta}$, for every $f, R(f) - R(f_{\theta}) = E_{X}[KL \Big(p_{\theta}(X) \parallel p_f(X) \Big)] \geq 0$ and therefore, $f_{\theta} = \operatorname*{arg\,min}_{f \in \operatorname{conv}(\cF) }R_{\pi_{\theta}}(f)$.
    
For any $A \subseteq \{1,\ldots,K\}$ let $p_{f,A}(x)=\sum_{k \in A} p_{f,k}(x)$.  
    Applying Lemma \ref{lemma:discrete-kl-lower-bound} we then have
        \[
    \begin{split}
       \inf & _{\hat{f}}\sup_{\pi}\mathbb{E}_{\pi^{\otimes n}}\Big[R(\hat{f})  - \inf_{f \in \operatorname{conv}(\cF)}R(f)\Big] \geq 
       \inf_{\hat{f}}\sup_{\pi_{\theta}:p_{\theta} \in P_{\tilde{\Theta}_{m}^{M}} }\mathbb{E}_{\pi^{\otimes n}}\Big[E_{X}\Big[KL\Big(p_{\theta}(X) \parallel  p_{\hat{f}}(X)\Big) \Big]
       \\
       & \geq
       2\inf_{\hat{f}}\sup_{\pi_{\theta}: p_{\theta} \in P_{\tilde{\Theta}_{m}^{M}} }\mathbb{E}_{\pi^{\otimes n}}
       \Big[\mathbb{E}_X\Big[(p_{\theta,A}(X)-p_{\hat{f},A}(X))^2\Big]\Big] 
        = 
        2\inf_{\hat{f}}\sup_{\pi_{\theta}: p_{\theta} \in P_{\tilde{\Theta}_{m}^{M}} }\mathbb{E}_{\pi^{\otimes n}}\Big[d^2(p_{\theta}, p_{\hat{f}}) \Big],
    \end{split}
    \]
    where the pseudo-metric  $d(p_{f_1},p_{f_2})=\sqrt{\mathbb{E}_X\Big[(p_{f_2,A}(X)-p_{f_1,A}(X))^2\Big]}$.

   Consider $A =\{1,3,5,\;\ldots,K-1\}$. Then, \eqref{equation:p_thetak} yields
        \begin{gather}\label{equation:conditional-distribution-shape-theorem-9}
            p_{\theta,A}(x) = \frac{e^{\epsilon \theta ^T x}}{e^{\epsilon \theta ^T x} + e^{-\epsilon \theta ^T x}}
    \end{gather}

    By a straightforward algebra, for any $s, t \in [-1,1]$, 
    \begin{gather}\label{equation:conditional-distribution-point-lower-bound-theorem-9}
        \frac{1}{2}|s-t|
        \geq
        \Big|\frac{e^s}{e^s+e^{-s}} - \frac{e^t}{e^t+e^{-t}}\Big|
               \geq
        \frac{1}{2\cosh^2(1)} |s-t|
    \end{gather}
    
Note that under the uniform distribution $\pi_X$ on $\Sigma$, $X_1,\ldots,X_M$ are i.i.d. Rademacher random variables.    
Combining \eqref{equation:hamming-distance-lower-bound-theorem-9}, \eqref{equation:conditional-distribution-shape-theorem-9} and \eqref{equation:conditional-distribution-point-lower-bound-theorem-9} for every $\theta, \tilde{\theta} \in \tilde{\Theta}^M_m$ we then have 
    \[
    \begin{split}
        d^2\Big(p_{\theta}, p_{\tilde{\theta}}\Big) 
        & =
        \mathbb{E}_{X} \Big[ (p_{\theta,A}(X) - p_{\tilde{\theta},A}(X))^2 \Big]
        \geq 
        \frac{1}{4\cosh^4(1)} \mathbb{E}_{X} \Big[(\epsilon \theta ^TX - \epsilon \tilde{\theta}^TX)^2 \Big]
        \\
        & =
        \frac{\epsilon^2}{4\cosh^4(1)} \mathbb{E}_{X}\Big[ \Big(\sum_{j=1}^M(\theta_j - \tilde{\theta}_j) X_j\Big)^2\Big]
         =\frac{\epsilon^2}{4\cosh^4(1)}\lVert \theta - \tilde{\theta} \rVert_2^2
         =
         \frac{\epsilon^2}{4\cosh^4(1)m^2}\rho(\theta, \tilde{\theta}) \\
         & \geq
         \frac{\epsilon^2}{16\cosh^4(1)m}
    \end{split}
    \]

    On the other hand, by \eqref{equation:p_thetak} and \eqref{equation:conditional-distribution-point-lower-bound-theorem-9}, for every $\theta, \tilde{\theta} \in \tilde{\Theta}^M_m$
        \begin{equation}\label{equation:upper-bound-kl-divergence-theorem-9}
    \begin{split}
        KL(\pi_{\theta} \parallel \pi_{\tilde{\theta}})
        &\leq
        \chi^2(\pi_{\theta}, \pi_{\tilde{\theta}})
        \\
        &=
        \sum_{x \in \Sigma}\sum_{k=1}^K\pi_{\tilde{\theta}}(x,k)
        \Big( \frac{\pi_{\theta}(x,k)}{\pi_{\tilde{\theta}}(x,k)}- 1 \Big)^2
        =
        \frac{1}{2^M}\sum_{x \in \Sigma}\sum_{k=1}^K
        \frac{\Big(p_{\tilde{\theta}, k}(x) - p_{\theta, k} (x)\Big)^2}
        {p_{\tilde{\theta}, k}(x)}
        \\
        &=
        \frac{1}{2^M}\sum_{x \in \Sigma}\Bigg[
        \frac{K}{2}
        \frac{\Big(\frac{2}{K}p_{\theta,A}(x)
        -\frac{2}{K}p_{\tilde{\theta},A}(x)\Big)^2}
        {\frac{2}{K}p_{\tilde{\theta},A}(x)}
        +
        \frac{K}{2}
        \frac{\Big(\frac{2}{K}p_{\theta,A}(-x)
        -\frac{2}{K}p_{\tilde{\theta},A}(-x)\Big)^2}
        {\frac{2}{K}p_{\tilde{\theta},A}(-x)}
        \Bigg]
        \\
        &\leq 
        \frac{1+e^2}{2^M}
        \sum_{x \in \Sigma}
        \Big[
        \big(p_{\theta,A}(x)-p_{\tilde{\theta},A}(x)\big)^2
        +
        \big(p_{\theta,A}(-x)-p_{\tilde{\theta},A}(-x)\big)^2
        \Big]
        \\
        &\leq
        \frac{9}{2^{M+1}}
        \sum_{x \in \Sigma}
        \Big(\epsilon\theta^Tx-\epsilon\tilde{\theta}^Tx\Big)^2
        =
        \frac{9}{2}
        E_X\Big[
        \big(\epsilon\theta^T X-\epsilon\tilde{\theta}^T X\big)^2
        \Big]
        \\
        &\leq
        \frac{9\epsilon^2}{m},
        \end{split}
    \end{equation}
    where we used the inequality $\frac{e^x}{e^x+e^{-x}} \geq \frac{1}{1+e^2}$ for $|x| \leq 1$.

    \noindent
    \textbf{Case 1.} $M \leq \sqrt{n}$. Take $m = \lfloor \frac{M}{2} \rfloor$ and $\epsilon^2 = \frac{c_1 \ln 3}{576}\frac{M^2}{n}$. Then, \eqref{equation:set-size-lower-bound-theorem-9} and \eqref{equation:upper-bound-kl-divergence-theorem-9} imply $KL(\pi_{\theta}^{\otimes n}\parallel \pi_{\tilde{\theta}}^{\otimes n}) 
        = 
        nKL(\pi_{\theta} \parallel \pi_{\tilde{\theta}}) 
        \leq \frac{1}{16}\ln|\tilde{\Theta}^M_m|$. 
    From Theorem 2.5 of \textcite{Introduction-to-Nonparametric-Estimation}, a variant of Fano lemma, it then follows that
    \begin{gather}\label{equation:lower-bound-convex-aggregation-small-M}
        \inf_{\hat{f}}\sup_{\pi}\mathbb{E}_{\pi^{\otimes n}}\Big[R(\hat{f}) - \inf_{f \in \operatorname{conv}(\cF)}R(f)  \Big]
        \geq
        c_2 \frac{M}{n}
    \end{gather}
    for some $c_2 > 0$.

    \noindent
    \textbf{Case 2.} $M > \sqrt{n}$. In this case
    take $m = \max\{1,\lfloor\frac{1}{2}\sqrt{\frac{n\ln 2}{\ln(1+\frac{M}{\sqrt{n}})}}\rfloor\}$ and $\epsilon^2 = \frac{c_1\ln 2}{2304}$. Applying the same arguments as in the previous case, Theorem 2.5 of \textcite{Introduction-to-Nonparametric-Estimation} implies now that
       \begin{gather}\label{equation:lower-bound-convex-aggregation-large-M}
        \inf_{\hat{f}}\sup_{\pi}\mathbb{E}_{\pi^{\otimes n}}\Big[R(\hat{f}) - \inf_{f \in \operatorname{conv}(\cF)}R(f) \Big]
        \geq
        c_3\sqrt{\frac{\ln(1 + \frac{M}{\sqrt{n}})}{n}}
    \end{gather}
    for some $c_3 > 0$

    Combining \eqref{equation:lower-bound-convex-aggregation-small-M} and \eqref{equation:lower-bound-convex-aggregation-large-M} completes the proof

\end{proof}

\begin{proof}[\textbf{Proof of Theorem \ref{theorem:continuous-CAEW-logistic-upper-bound}}]

For \(i=0,\ldots,n\), define
\[
C_i=
\int_{\Delta_M}e^{-iL^i(f_\theta)}\,d\theta,
\]
where recall that \(d\theta\) denotes the \((M-1)\)-dimensional volume measure on \(\Delta_M\), and \(C_0=\int_{\Delta_M}d\theta\). Let
\[
w_i(\theta)=\frac{e^{-iL^i(f_\theta)}}{C_i}
\]

For \(i=0,\ldots,n-1\), consider the probability mixture
\[
q_i(x)=
\int_{\Delta_M}
\operatorname{softmax}(f_\theta(x))
w_i(\theta)\,d\theta.
\]

Since \(w_i\) is a probability density on \(\Delta_M\),  \(q_i(x)\in\Delta_K\).
By the definition of \(\hat w_{CAEW}\),
$$
\operatorname{softmax}(\hat f_{CAEW}^{cP}(x))
=
\int_{\Delta_M}
\operatorname{softmax}(f_\theta(x))
\hat w_{CAEW}(\theta)\,d\theta
=
\frac1n\sum_{i=0}^{n-1}q_i(x)
$$

The first equality holds because the vector inside the componentwise logarithm is a probability vector whose components sum to one.
By convexity of \(-\ln(\cdot)\),
\begin{equation} \label{equation:A}
R(\hat f_{CAEW}^{cP})
=
\mathbb E_\pi\left[
-\ln\left(
\frac1n\sum_{i=0}^{n-1}q_{i,Y}(X)
\right)\right]
\le
\frac1n\sum_{i=0}^{n-1}
R(\ln q_i)
\end{equation}
and therefore,
\[
\mathbb E_{\pi^{\otimes n}}
R(\hat f_{CAEW}^{cP})
\le
\frac1n\sum_{i=0}^{n-1}
\mathbb E_{\pi^{\otimes n}}R(\ln q_i).
\]

Since \(q_i\) depends only on the first \(i\) observations and \((X_{i+1},Y_{i+1})\) is an independent observation from \(\pi\), the definition of risk and the law of iterated expectations imply
\[
\mathbb E_{\pi^{\otimes n}}R(\ln q_i)
=
\mathbb E_{\pi^{\otimes n}}
[-\ln q_{i,Y_{i+1}}(X_{i+1})]
\]

Moreover,
\[
q_{i,Y_{i+1}}(X_{i+1})=
\frac1{C_i}
\int_{\Delta_M}
e^{-iL^i(f_\theta)}
\operatorname{softmax}(f_\theta(X_{i+1}))_{Y_{i+1}}
\,d\theta
=\frac1{C_i}
\int_{\Delta_M}
e^{-(i+1)L^{i+1}(f_\theta)}\,d\theta
=\frac{C_{i+1}}{C_i},
\]
where we used
$\operatorname{softmax}(f_\theta(x))_y
=e^{-l(f_\theta(x),y)}$
and additivity of the loss:
$(i+1)L^{i+1}(f_\theta)
=
iL^i(f_\theta)+l(f_\theta(X_{i+1}),Y_{i+1})$.

Consequently,
\[
\frac1n\sum_{i=0}^{n-1}
-\ln q_{i,Y_{i+1}}(X_{i+1})
=
-\frac1n\sum_{i=0}^{n-1}
\ln\Bigl(\frac{C_{i+1}}{C_i}\Bigr)
=-\frac1n\ln\Bigl(\frac{C_n}{C_0}\Bigr)
\]

Combining this identity with \eqref{equation:A} yields
\begin{equation} \label{equation:CAEW-sequential}
\mathbb E_{\pi^{\otimes n}}R(\hat f_{CAEW}^{cP})
\le
-\frac1n\mathbb E_{\pi^{\otimes n}}
\ln\Bigl(\frac{C_n}{C_0}\Bigr)
\end{equation}
\ignore{
where recall that
\[
C_n=\int_{\Delta_M}e^{-nL_n(f_\theta)}\,d\theta,
\qquad
C_0=\int_{\Delta_M}d\theta
\]
}

Let \(\theta^*\) be a minimizer of \(R(f_\theta)\) over \(\Delta_M\). For any \(0 <\varepsilon\le1\), define
\[
\Delta_M^\varepsilon
=(1-\varepsilon)\theta^*+\varepsilon\Delta_M
\]

Since \(\Delta_M^\varepsilon\subseteq\Delta_M\),
\[
C_n\ge
\int_{\Delta_M^\varepsilon}
e^{-nL^n(f_\theta)}\,d\theta
\]
and by Jensen's inequality,
\[
C_n
\ge
\operatorname{Vol}(\Delta_M^\varepsilon)
\exp\left(
-\frac{n}{\operatorname{Vol}(\Delta_M^\varepsilon)}
\int_{\Delta_M^\varepsilon}
L^n(f_\theta)\,d\theta
\right)
\]
Since \(\Delta_M\) is an \((M-1)\)-dimensional simplex, the affine contraction by a factor \(\varepsilon\) multiplies its volume by \(\varepsilon^{M-1}\), that is,
\[
\operatorname{Vol}(\Delta_M^\varepsilon)
=\varepsilon^{M-1}C_0,
\]
implying
\[
-\frac1n\ln\Bigl(\frac{C_n}{C_0}\Bigr)
\le
\frac{1}{\operatorname{Vol}(\Delta_M^\varepsilon)}
\int_{\Delta_M^\varepsilon}L^n(f_\theta)\,d\theta
+
\frac{M-1}{n}\ln\frac1\varepsilon
\]

Taking expectations with respect to the training sample gives
\begin{equation} \label{equation:CAEW-Jensen}
-\frac1n\mathbb E_{\pi^{\otimes n}}
\ln\Bigl(\frac{C_n}{C_0}\Bigr)
\quad\le
\frac{1}{\operatorname{Vol}(\Delta_M^\varepsilon)}
\int_{\Delta_M^\varepsilon}R(f_\theta)\,d\theta
+
\frac{M-1}{n}\ln\frac1\varepsilon
\end{equation}

Every \(\theta\in\Delta_M^\varepsilon\) can be written as 
\(\theta=(1-\varepsilon)\theta^*+\varepsilon\vartheta\), where
$\vartheta\in\Delta_M$. Hence, by definition of $f_\theta$,
\[
f_\theta=(1-\varepsilon)f_{\theta^*}+\varepsilon f_\vartheta,
\]
and by convexity of the logistic loss,
\[
R(f_\theta)
\le (1-\varepsilon)R(f_{\theta^*})
+\varepsilon R(f_\vartheta)
=R(f_{\theta^*})
+\varepsilon\bigl[R(f_\vartheta)-R(f_{\theta^*})\bigr]
\]

Since all score functions in \(\operatorname{conv}(\mathcal F)\) are bounded by \(B\) in sup-norm, 
their logistic risks differ by at most \(4B\)
(see \eqref{equation:bounded_loss}). Thus,
\[
R(f_\theta)\le R(f_{\theta^*})+4B\varepsilon
\]

Combining the preceding inequality with \eqref{equation:CAEW-sequential} and \eqref{equation:CAEW-Jensen} yields
\begin{equation} \label{equation:F}
\mathbb E_{\pi^{\otimes n}}R(\hat f_{CAEW}^{cP})
-\inf_{f\in\operatorname{conv}(\cF)}R(f)
\le
4B\varepsilon
+\frac{M-1}{n}\ln\frac1\varepsilon
\end{equation}

Choose
\[
\varepsilon=
\frac{M-1}{M-1+4Bn}.
\]

Then
\[
4B\varepsilon
=
\frac{4B(M-1)}{M-1+4Bn}
\le\frac{M-1}{n}
\]
and
\[
\ln\frac1\varepsilon
=
\ln\left(1+\frac{4Bn}{M-1}\right).
\]

Substituting into \eqref{equation:F} yields
\[
\mathbb E_{\pi^{\otimes n}}
[R(\hat f_{CAEW}^{cP})]
\le
\inf_{f\in\operatorname{conv}(\cF)}R(f)
+
\frac{M-1}{n}
\left[
1+\ln\left(1+\frac{4Bn}{M-1}\right)
\right]
\]
This completes the proof.

\end{proof}

\begin{proof}[\textbf{Proof of Theorem \ref{theorem:probability-mixture-convex-aggregation-erm}}]
The proof follows arguments similar to those used for Theorem \ref{theorem:erm-convex-hull-upper-bound}.

Define a vector $P_y(x)=\left(p_{1y}(x),\ldots,p_{My}(x)\right)^T \in \mathbb{R}^M$ and a matrix $W_y(x)=P_y(x)P_y(x)^T \in \mathbb{R}^{M \times M}$.  We have
    \begin{gather*}
         \nabla_\theta l_{prob}(p_\theta(x),y)= -\frac{1}{\theta^TP_y(x)}P_y(x) 
         \\
         \nabla_\theta^2(l_{prob}(p_\theta(x),y))=\frac{1}{(\theta^TP_y(x))^2}W_y(x)
     \end{gather*}

For any $\theta_1,\theta_2\in\Delta_M$, Taylor's expansion yields
\[
l_{\mathrm{prob}}(p_{\theta_1}(x),y)
-l_{\mathrm{prob}}(p_{\theta_2}(x),y)
=
\nabla_\theta l_{\mathrm{prob}}(p_{\theta_2}(x),y)^\top(\theta_1-\theta_2)
+\frac{1}{2(\tilde\theta^\top P_y(x))^2}
||\theta_1-\theta_2||_{W_y(x)}^2,
\]
where $\tilde\theta$ lies on the segment joining $\theta_1$ and $\theta_2$.

Since $e^{-2B}/K\le\theta^\top P_y(x)\le e^{2B}/K$ for every $\theta\in\Delta_M$, it follows that
\[
l_{\mathrm{prob}}(p_{\theta_1}(x),y)
-l_{\mathrm{prob}}(p_{\theta_2}(x),y)
\ge
\nabla_\theta l_{\mathrm{prob}}(p_{\theta_2}(x),y)^\top(\theta_1-\theta_2)
+\frac{K^2e^{-4B}}2||\theta_1-\theta_2||_{W_y(x)}^2,
\]
\ignore{
and
\[
l{\mathrm{prob}}(p_{\theta_1}(x),y)
-l_{\mathrm{prob}}(p_{\theta_2}(x),y)\le
\nabla_\theta l_{\mathrm{prob}}(p_{\theta_2}(x),y)^\top(\theta_1-\theta_2)
+\frac{K^2e^{4B}}2|\theta_1-\theta_2|_{W_y(x)}^2.
\]
}

On the other hand, by the mean value theorem, for some $\bar\theta$ on the segment joining $\theta_1$ and $\theta_2$,
\[
l_{prob}(p_{\theta_1}(x),y)-l_{prob}(p_{\theta_2}(x),y)
=-\frac{(\theta_1-\theta_2)^TP_y(x)}{\bar\theta^TP_y(x)}.
\]
Consequently,
\[
\Bigl(l_{prob}(p_{\theta_1}(x),y)-l_{prob}(p_{\theta_2}(x),y)\Bigr)^2
\le K^2e^{4B}|\theta_1-\theta_2|_{W_y(x)}^2.
\]

Thus, the multinomial logistic loss in probability space is $K^2e^{-4B}$-strongly convex and $Ke^{2B}$-Lipschitz over $\Delta_M$ with respect to the seminorm induced by $W_y(x)$.

\ignore{
Then, since $\frac{e^{-2B}}{K} \leq p_{jy}(x) \leq \frac{e^{2B}}{K}$ for all $j$ and $y$,
$$    
\Bigl(l_{prob}(p_{\theta_1}(x),y)-l_{prob}(p_{\theta_2}(x),y)\Bigr)^2
        =\frac{1}{(\tilde{\theta}^TP_y(x))^2}\lVert \theta_1-\theta_2 \rVert^2_{W_y(x)} \leq K^2 e^{4B} \lVert \theta_1-\theta_2 \rVert^2_{W_y(x)}
$$
for some $\tilde{\theta} \in \Delta_M$.

Similarly, by Taylor's expansion
\begin{equation*}
\begin{split}
l_{prob}(p_{\theta_1}(x),y) - l_{prob}(p_{\theta_2}(x),y) &  =
         \nabla_\theta l_{prob}(p_{\theta_2}(x),y)^T(\theta_1-\theta_2)
         +
         \frac{1}{2}(\theta_1-\theta_2)^T \frac{1}{(\hat{\theta}^TP_y(x))^2}W_y(x)(\theta_1-\theta_2) 
         \\
         & \geq 
         \nabla_\theta l_{prob}(p_{\theta_2}(x),y)^T(\theta_1-\theta_2)
         +
         \frac{K^2e^{-4B}}{2} \lVert \theta_1-\theta_2 \rVert_{W_y(x)}^2
\end{split}
\end{equation*}
for some $\hat{\theta} \in \Delta_M$.
}

Therefore, for the empirical risk $L^n_{prob}(p_\theta)=\frac{1}{n}\sum_{i=1}^n l_{prob}(p_\theta(x_i),y_i)$ we have
\begin{gather*}
    \frac{1}{n}\sum_{i=1}^n \Big( l_{prob}(p_{\theta_1}(x_i),y_i) - l_{prob}(p_{\theta_2}(x_i),y_i) \Big)^2
    \leq 
        K^2 e^{4B} \lVert \theta_1-\theta_2 \rVert_W^2 \\
 L^n_{prob}(p_{\theta_1}) - L^n_{prob}(p_{\theta_2})   \geq 
        \nabla_\theta L^n_{prob}(p_{\theta_2})^T(\theta_1-\theta_2)
         +
         \frac{K^2e^{-4B}}{2} \lVert \theta_1-\theta_2 \rVert_W^2,
\end{gather*}
where the matrix $W=\frac{1}{n} \sum_{i=1}^n W_{y_i}(x_i)$.

Hence, $L^n_{prob}(p_\theta)$ satisfies Assumption 3.1 of \textcite{puchkin2023exploringlocalnormsexpconcave} with $L^2=K^2e^{4B}$ and $\sigma = K^2 e^{-4B}$, and the seminorm induced by $W$. Applying Theorem 3.4 in the same paper then implies that with probability at least $1-\delta$
     \[
        R_{prob}(\hat{p}^{l_{prob},\operatorname{conv}(\cP)}_{ERM}) - \min_{p \in \operatorname{conv}(\cP)}R_{prob}(p) \leq c_1\frac{e^{8B}(M + \ln(\frac{1}{\delta}))}{n}
     \]   
for any $0 < \delta < 1$. Integrating the tail probability yields
 \begin{gather}\label{equation:probability-mixture-erm-upper-bound-small-M}
        \mathbb{E}_{\pi^{\otimes n}}\Big[R_{prob}(\hat{p}^{l_{prob},\operatorname{conv}(\cP)}_{ERM})\Big] 
        - 
        \min_{p \in \operatorname{conv}(\cP)}R_{prob}(p) 
        \leq 
        c_1(B)\frac{M}{n}
     \end{gather}

To establish the $\sqrt{\ln(M)/n}$ bound, define $g_j(x,y)=Kp_{jy}(x)$, $j=1,\ldots,M$. Since $e^{-2B}\le g_j(x,y)\le e^{2B}$, we have
\[
l_{prob}(p_\theta(x),y)
=\ln K-\ln\left(\sum_{j=1}^M\theta_jg_j(x,y)\right).
\]
The additive constant $\ln K$ does not affect the excess risk, while $-\ln(\cdot)$ is $e^{2B}$-Lipschitz on $[e^{-2B},e^{2B}]$. Thus, by the standard symmetrization inequality, the contraction principle, and Massart's lemma,
\[
\mathbb E_{\pi^{\otimes n}}
\left[R_{prob}(\hat p_{ERM}^{l_{prob},\operatorname{conv}(\cP)})\right]
-\min_{p\in\operatorname{conv}(\cP)}R_{prob}(p)
\le c_2(B)\sqrt{\frac{\ln M}{n}}.
\]

\ignore{
The proof of the $\sqrt{\frac{\ln(M)}{n}}$-term
follows the same symmetrization–contraction–Massart arguments as in Theorem \ref{theorem:erm-convex-hull-upper-bound} with the adjusted Lipschitz constant $Ke^{2B}$ for $-\ln(\cdot)$ on $[\frac{e^{-2B}}{K},\frac{e^{2B}}{K}]$
}
    
\end{proof}

\begin{proof}[\textbf{Proof of Theorem~\ref{theorem:probability-mixture-convex-aggregation-lower-bound}}]
The proof goes along the lines of the proof of Theorem \ref{theorem:convex-aggregation-lower-bound}.

    Assume for simplicity that $K$ is even. Consider $x \in \Sigma=\{-1,1\}^M$ and assume the uniform prior $\pi_X(x)=\frac{1}{2^M},\; x\in\Sigma$.

Define the vectors
$u=\bigl(\tfrac{1}{K},\ldots,\tfrac{1}{K}\bigr)^\top\in\mathbb{R}^K$
and
$v=\bigl(
\underbrace{\tfrac{1}{K},\ldots,\tfrac{1}{K}}_{\scriptscriptstyle K/2},
\underbrace{-\tfrac{1}{K},\ldots,-\tfrac{1}{K}}_{\scriptscriptstyle K/2}
\bigr)^\top\in\mathbb{R}^K$,
and the functions
$g_j(x)=x_j v$, $j=1,\ldots,M$,
where $x \in \Sigma$.

Define the dictionary of probability vectors
\[
p_j(x)=u+\epsilon g_j(x),\qquad j=1,\ldots,M,
\]
where $0<\epsilon \leq \frac{1}{2}$ will be defined later (note that since $\sum_{k=1}^K v_k=0$ and $|v_k|=\frac{1}{K}$, 
$p_j(x)\in\Delta_K$ for every $x$ and every $j$). The corresponding dictionary $\cF$ consists of score functions $f_j(x) \in \mathbb{R}^K,\;j=1,\ldots,M$, with 
$f_{jk}(x)=\ln p_{jk}(x)-\overline{\ln p_{j\cdot}(x)}$~.

For every $\theta\in \Delta_M$ we have
    \begin{gather}\label{equation:p-theta-calculation}
        p_\theta(x)
        =
        \sum_{j=1}^M\theta_j p_j(x) 
        = 
        \sum_{j=1}^M\theta_j(u + \epsilon g_j(x))
               =
        u + \epsilon v \sum_{j=1}^M\theta_j x_j
    \end{gather}

    Take $A = \{1,\ldots, \frac{K}{2}\}$. Then,
    \begin{gather} \label{equation:p_thetaA}
        p_{\theta,A}(x) = \sum_{k \in A} p_{\theta,k}(x)=\frac{1 + \epsilon \sum_{j=1}^M\theta_j x_j}{2}
    \end{gather}

      As we argued in the proof of Theorem \ref{theorem:convex-aggregation-lower-bound}, there exists a subset of $m$-sparse vectors $\tilde{\Theta}^M_m \subseteq  \{\theta \in \{0,\frac{1}{m}\}^M :   \lVert \theta \rVert_0 = m\}$ such that for any $\theta_1, \theta_2 \in\tilde{\Theta}^M_m$
    \begin{gather}\label{equation:lower-bound-hamming-distance}
        \rho(\theta_1, \theta_2) \geq \frac{m}{4}
    \end{gather}
    \begin{gather}\label{equation:lower-bound-set-size}
        \ln|\tilde{\Theta}^M_m| \geq c_1m\ln(1+M/m)
    \end{gather}

Consider the distribution $\pi_{\theta}(x,k)=\pi_{X}(x) p_{\theta,k}(x)$ on $\Sigma \times \mathcal{Y}$. Following the arguments in the proof of Theorem \ref{theorem:convex-aggregation-lower-bound},
    since $\tilde{\Theta}_m^M \subset \Delta_M$, by Lemma \ref{lemma:discrete-kl-lower-bound},
    \begin{gather*}
        \inf_{\hat{{p}}}\sup_{\pi}\mathbb{E}_{\pi^{\otimes n}}\Big[R_{prob}(\hat{p}) - \inf_{p \in \operatorname{conv}(\cP)}R_{prob}(p)\Big] 
        \geq
        \ignore{
            \inf_{\hat{{p}}}\sup_{\pi: p_\theta \in p_{\tilde{\Theta}^M_m}}\mathbb{E}_{\pi^{\otimes n}} \Big[\mathbb{E}_X[KL(\pi(x) \parallel \hat{p}(x))]\Big] \\
            \geq 
            2\inf_{\hat{{p}}}\sup_{\pi: \pi_{x} \in p_{\tilde{\Theta}^M_m}}\mathbb{E}_{\pi^{\otimes n}}\Big[\mathbb{E}_X\Big[(\pi_A(x)-\hat{p}_A(x))^2\Big]\Big] 
            = 
        }
        2\inf_{\hat{{p}}(x)}\sup_{\pi_\theta: p_\theta \in p_{\tilde{\Theta}^M_m}}\mathbb{E}_{\pi^{\otimes n}}\Big[d^2(p_\theta, \hat{p}) \Big],
    \end{gather*}
    where the pseudo-metric $d(p_\theta,\hat{p})=\sqrt{\mathbb{E}_X\Big[(\hat{p}_A(X)-p_{\theta,A}(X))^2\Big]}$. 

Note that, under the uniform $\pi_X$ on $\Sigma$, $X_1,\ldots,X_M$  are i.i.d. Rademacher random variables.
   Therefore, \eqref{equation:p_thetaA}-\eqref{equation:lower-bound-hamming-distance} imply that for every $\theta, \tilde{\theta} \in \tilde{\Theta}^M_m$, 
    \begin{gather*}\label{equation:distnace-as-expected-value}
        d^2\Big(p_{\theta}, p_{\tilde{\theta}}\Big) 
        =
        \mathbb{E}_{X} \Big[ (p_{\theta,A}(X) - p_{\tilde{\theta},A}(X))^2 \Big]
        =
        \frac{\epsilon^2}{4}\mathbb{E}_X \Big[ \Big(\sum_{j=1}^M(\theta_j - \tilde{\theta}_j) X_j \Big)^2 \Big]
        =
        \frac{\epsilon^2}{4} \lVert \theta - \tilde{\theta} \rVert_2^2 \\
         =
        \frac{\epsilon^2}{4m^2}\rho(\theta, \tilde{\theta})
        \geq 
        \frac{\epsilon^2}{16m}
    \end{gather*}
        
    On the other hand, by \eqref{equation:p-theta-calculation},
    \begin{equation}\label{equation:upper-bound-kl-divergence}
    \begin{split}
        KL(\pi_{\theta} \parallel \pi_{\tilde{\theta}}) 
        & \leq 
        \chi^2(\pi_{\theta}, \pi_{\tilde{\theta}})
        =
        \sum_{x \in \Sigma}\sum_{k=1}^K\pi_{\tilde{\theta}}(x,k)\Big( \frac{\pi_{\theta}(x,k)}{\pi_{\tilde{\theta}}(x,k)}- 1 \Big)^2
        = 
        \frac{1}{2^M}\sum_{x \in \Sigma}\sum_{k=1}^K\frac{\Big(p_{\tilde{\theta},k}(x)  -  p_{\theta,k}(x)\Big)^2}{p_{\tilde{\theta},k}(x)} \\ 
        & = \frac{\epsilon^2}{2^{M}}\sum_{x \in \Sigma}\frac{\Big(\sum_{j}(\tilde{\theta}_j- \theta_j)x_j\Big)^2}{1-\epsilon^2\Big( \sum_{j}\tilde{\theta}_jx_j \Big)^2} 
        \leq 
        2\epsilon ^ 2\mathbb{E}_X \Big[ \Big(\sum_{j}(\tilde{\theta}_j- \theta_j)X_j\Big)^2 \Big]
             =
             2\epsilon^2 \lVert \tilde{\theta}-\theta \rVert_2^2 \\ 
            & \leq 
             \frac{4\epsilon^2}{m}
    \end{split}
    \end{equation}

    \noindent
    \textbf{Case 1.} $M \leq \sqrt{n}$. Take $m = \lfloor\frac{M}{2} \rfloor$ and $\epsilon^2 = \frac{c_1 \ln 3}{512}\frac{M^2}{n}$. Then, 
    \eqref{equation:lower-bound-set-size}-\eqref{equation:upper-bound-kl-divergence} imply $KL(\pi_{\theta}^{\otimes n} \parallel \pi_{\tilde{\theta}}^{\otimes n}) 
        = 
        nKL(\pi_{\theta} \parallel \pi_{\tilde{\theta}}) 
        \leq \frac{1}{16}\ln|\tilde{\Theta}^M_m|$
    and by Theorem 2.5 of \textcite{Introduction-to-Nonparametric-Estimation},
    \begin{gather}\label{equation:lower-bound-probability-mixture-convergence-small-M}
        \inf_{\hat{{p}}}\sup_{\pi}\mathbb{E}_{\pi^{\otimes n}}\Big[R_{prob}(\hat{p}) - \inf_{p \in \operatorname{conv}(\cP)}R_{prob}(p) \Big]
        \geq
        c_2\frac{M}{n}
    \end{gather}
    for some $c_2 > 0$.
    
    \noindent
    \textbf{Case 2.} $M > \sqrt{n}$. In this case
    take $m = \max\{1,\lfloor\frac{1}{2}\sqrt{\frac{n\ln 2}{\ln(1+\frac{M}{\sqrt{n}})}}\rfloor\}$ and $\epsilon^2 = \min \{\frac{c_1 \ln 2}{1024}, \frac{1}{4} \}$. Applying the same arguments as in the previous case, Theorem 2.5 of \textcite{Introduction-to-Nonparametric-Estimation} implies now that
       \begin{gather}\label{equation:lower-bound-probability-mixture-convergence-large-M}
        \inf_{\hat{{p}}}\sup_{\pi}\mathbb{E}_{\pi^{\otimes n}}\Big[R_{prob}(\hat{p}) - \inf_{p \in \operatorname{conv}(\cP)}R_{prob}(p) \Big]
        \geq
        c_3\sqrt{\frac{\ln(1 + \frac{M}{\sqrt{n}})}{n}}
    \end{gather}
    for some $c_3 > 0$.

Combining \eqref{equation:lower-bound-probability-mixture-convergence-small-M} and \eqref{equation:lower-bound-probability-mixture-convergence-large-M} completes the proof.

\end{proof}

\printbibliography

\end{document}